\documentclass[matann2]{svjour}

\usepackage{amsmath}
\usepackage{amsfonts}
\usepackage{latexsym}
\usepackage{mathrsfs}
\usepackage{times}
\usepackage{array}
\usepackage{amscd}
\usepackage{a4}
\usepackage[mathcal]{euscript}
\usepackage{amssymb}
\usepackage{pstricks}

\input{boxedeps} 
\SetEPSFDirectory{} 
\HideDisplacementBoxes
\SetRokickiEPSFSpecial  
\DeclareMathAlphabet{\ams}{U}{msb}{m}{n}
\DeclareMathAlphabet{\goth}{U}{euf}{m}{n}

\def\sl{\text{SL}}

\def\gl{\text{GL}}

\def\m{\text{M}}
\def\d{\text{D}}

\def\sp{\text{Sp}}

\def\rk{\text{rk}\,}

\def\im{\text{im}\,}
\def\aut{\text{Aut}}

\def\ov{\overline}

\def\wh{\widehat}

\def\dom{\text{dom}}
\def\codim{\text{codim}\,}

\def\II{\mathscr I}
\def\JJ{\mathscr J}

\def\BB{\mathcal B}
\def\CC{\mathcal C}
\def\FF{\mathcal F}
\def\AA{\mathcal A}

\def\HH{\mathcal H}

\def\PP{\mathcal P}
\def\QQ{\mathcal Q}
\def\TT{\mathcal T}

\def\SSS{\mathcal S}

\def\BBB{\goth{B}}
\def\XXX{\goth{X}}

\def\ve{\varepsilon}

\def\aa{\alpha}
\def\ww{\omega}

\def\vphi{\varphi}

\def\ve{\varepsilon}
\def\Om{\Omega}

\def\wh{\widehat}

\def\Z{\ams{Z}}
\def\R{\ams{R}}
\def\C{\ams{C}}\def\Q{\ams{Q}}
\def\F{\ams{F}}
\def\T{\ams{T}}

\def\G{\ams{G}}
\def\M{\ams{M}}

\def\vv{\mathbf{v}}
\def\uu{\mathbf{u}}
\def\xx{\mathbf{x}}

\def\quo{/\kern -.45em\sim}

\newpsobject{showgrid}{psgrid}{subgriddiv=1,griddots=5,gridlabels=6pt}

\def\ds{\displaystyle}

\def\Langle{\langle\kern -2pt\langle}
\def\Rangle{\rangle\kern -1.9pt\rangle}
\title{Partial mirror symmetry I: reflection monoids}

\author{Brent Everitt and John Fountain
\thanks{The first author would like to thank 
Bob Howlett 
for helpful suggestions.
Some of the results of this paper were 
obtained while he was visiting
the Institute for Geometry and its Applications, 
University of Adelaide, Australia. 
He is grateful for their hospitality. The second author would like to
thank Mark Kambites for some helpful discussions, and Chris Hollings for
his help with several calculations. A grant from the Royal Society made
it possible for him to visit the University of Adelaide to continue the
work reported here. He would like to express his gratitude to the members of
the Glenelg Mathematics Institute for their kindness and hospitality
during his visit to Adelaide.}
}

\institute{
Department of Mathematics, University of York, York
YO10 5DD, United Kingdom. \email{bje1@york.ac.uk} (Brent Everitt), 
\email{jbf1@york.ac.uk} (John Fountain).
}

\titlerunning{}
\authorrunning{Brent Everitt and John Fountain}

\begin{document}

\maketitle


\begin{abstract}
This is the first of a series of papers in which 
we initiate and develop the theory of reflection monoids, 
motivated by the theory of reflection
groups. The main results 
identify a number of important inverse semigroups as reflection monoids,
introduce new examples, and determine their orders.
\end{abstract}


\section*{Introduction}

The symmetric group $\goth{S}_n$ comes in many guises: as the permutation group of the
set $\{1,\ldots,n\}$; as the group
generated by reflections in the hyperplanes $x_i-x_j=0$ of 
an $n$-dimensional Euclidean space;
as the Weyl group 
of the reductive algebraic group $\gl_n$, or (semi)simple group $\sl_{n+1}$,
or simple Lie algebra $\goth{sl}_{n+1}$; as the Coxeter group
associated to Artin's braid group, and so on.

If one thinks of $\goth{S}_X$ as the group of (global) symmetries of $X$, then the partial
symmetries naturally lead one to consider the \emph{symmetric inverse monoid} 
$\II_X$, whose
elements are the partial bijections $Y\rightarrow Y'$ $(Y,Y'\subset X)$.
It too has many other faces. It arises in its incarnation 
as the ``rook monoid'' as 
the so-called Renner monoid of the reductive algebraic monoid $\m_n$
(see \S \ref{subsection:rennermonoids} for the definitions). An 
associated Iwahori theory and representations 
have been worked out by Solomon \cite{Solomon02,Solomon90}. 
There is a braid connection too, with $\II_n$ naturally associated to 
the inverse monoid of ``partial braids'' defined recently in \cite{Easdown04}.

But what is missing from all this is a realization of $\II_n$ as some kind
of ``partial'' reflection monoid, or indeed, a definition and theory
of partial mirror symmetry and the monoids generated by partial reflections
that generalizes the theory of reflection groups.

Such is the purpose of the present paper. Reflection monoids are defined
as monoids generated by certain partial linear isomorphisms 
$\aa:X\rightarrow Y$ ($X,Y$ subspaces of $V$), 
that are the
restrictions (to $X$) of reflections. Initially one is faced with many 
possibilities, with the challenge being to impose enough structure for
a workable theory while still encompassing as many interesting examples as 
possible. 
It turns out that a solution is to consider monoids of partial linear
isomorphisms whose domains form a $W$-invariant semilattice for
some reflection group $W$ acting on $V$. 

Two pieces of data will characterise a reflection monoid:
a reflection \emph{group\/} and a collection of well behaved domain subspaces
(see \S \ref{section:reflectionmonoids} for the precise definitions).
What results is a theory of reflection monoids for which our main
theorems determine their orders, 
presentations and identify the natural examples
(it turns out that much of the general theory works when an arbitrary group is fed
into the input data, but at various crucial stages the reflection group 
structure will be used in a essential way to obtain results in specific examples).

For instance, just as $\goth{S}_n$ is the reflection group associated to the 
type $A$ root system, so
now $\II_n$ becomes the reflection monoid associated to the type $A$ root
system, and where the domains
form a Boolean lattice (see \S\ref{subsection:permutationmonoids}). The analogy
continues: the group of signed permutations
of $\{1,\ldots,n\}$ is the Weyl group of type $B$, and the inverse monoid
$\JJ_n$ of \emph{partial\/} signed permutations becomes the reflection 
monoid of type $B$,
with again this Boolean lattice of domain subspaces. 

By the ``rigidity of tori'', a maximal torus $T$ in a linear algebraic group
$\G$ has automorphisms a finite group, the Weyl group of $\G$, 
and this is a reflection
group in the space $\XXX(T)\otimes\R$, where $\XXX(T)$ is the character group
of the torus. A similar role is played in the theory of linear algebraic
\emph{monoids\/} by the Renner monoid (see \S\ref{subsection:rennermonoids} for 
the definitions). 
One might hope that the Renner monoids are examples of reflection monoids, but
in fact it turns out to be more complicated than this. We construct a reflection monoid
in $\XXX(T)\otimes\R$, where the extra piece of data, the semilattice of domain spaces, 
comes from the character semigroup $\XXX(\ov{T})$ of the Zariski
closure of $T$. This reflection monoid then maps homomorphically onto the
Renner monoid, with the two isomorphic in some cases.

Another interesting class of examples arises from the theory of hyperplane arrangements.
The 
\emph{reflection arrangement monoids\/} have as their input data a reflection group
and for the domains, the intersection lattice of the reflecting hyperplanes.
These intersection lattices possess many beautitful combinatorial and 
algebraic properties (see \cite{Orlik92}).
Thus, the reflection arrangement monoids tie up reflection groups and the intersection
lattices of their reflecting hyperplanes in one very natural algebraic object.

This first paper has been written so as to include in its readership
workers in both reflection groups and semigroups, and
is organized as follows: \S \ref{section:reflectiongroups}
contains background material on reflection groups; \S 
\ref{section:systems} introduces the semilattice of subspaces
forming the domains of our partial isomorphisms, and discusses in some detail
two classes of examples arising from hyperplane arrangements. Reflection monoids
proper are defined in \S \ref{section:reflectionmonoids}, along with
basic concepts in semigroup theory, and a number of their
basic properties are considered. The final section gives
three families
of examples:
the Boolean, Renner
and reflection arrangement monoids along with their orders in a number
of cases.

In the sequel \cite{Everitt:Fountain} to this paper,
a general presentation is derived (among other things) using 
the factorizable inverse monoid structure, and interpreting the various ingredients
of a presentation for such given recently in \cite{Easdown05}.
This presentation is determined explicity (and massaged a little more)
for the Boolean and arrangement monoids associated to 
the classical Weyl groups. The benchmark here is provided by a classical 
presentation \cite{Popova61} for the symmetric inverse monoid $\II_n$, 
which we rederive
in its new guise as the ``Boolean monoid of type $A$''.

\section{Preliminaries from reflection groups}\label{section:reflectiongroups}

Before venturing into partial mirror symmetry, we summarize the results we will
need from (full) mirror symmetry, ie: from the theory of reflection groups.
A number of these will not be needed until the sequel 
\cite{Everitt:Fountain} to this paper, but we place them here for convenience.
Standard references are \cite{Bourbaki02,Humphreys90}, and more recently
\cite{Kane01}.

Let $\F$ be a field, $V$ an $\F$-vector space and $GL(V)$ the group of
linear isomorphisms $V\rightarrow V$. A {\em reflection\/} is a 
non-trivial element
of finite order in $GL(V)$ that is 
semisimple and leaves pointwise invariant a hyperplane $H\subset V$.
A subgroup $W\subset GL(V)$ is a {\em reflection
group\/} when it is generated by reflections.

The most commonly studied 
examples arise in the cases $\F=\R,\C,\F_q$ and $\Q_p$ ($p$-adics),
and as all but one of the eigenvalues of an order $n$ reflection are equal to $1$, 
the last must be a primitive $n$-th root of unity in $\F$.
Thus $\F$ plays a role in the kinds of orders that 
reflections may have: they are involutions in the reals and 
$2$-adics, can have arbitrary finite order in the complexes, order dividing
$p-1$ in $\Q_p$ for $p$ an odd prime, and so on.

There are classical and celebrated classifications due to Coxeter 
\cite{Coxeter35,Coxeter34} in the reals, Shephard-Todd
\cite{Shephard54} (complexes), Clark-Ewing \cite{Clark74} ($p$-adics)
and Wagner \cite{Wagner81,Wagner80}, 
Zalesski{\u\i}-Sere{\v{z}}kin \cite{Zalesskii80} ($\F_q$).
In this paper, more for concreteness than any other reason, we will restrict
ourselves to $\F=\R$ and $\C$, and to reflection groups $W$ that are finite. 


Any finite subgroup of $GL(V)$ for $V$ a complex space leaves invariant
a positive definite Hermitian form, obtained in the usual way by an averaging
process. 
Two such reflection groups $W_i\subset GL(V_i)$ are
isomorphic if and only if there is a vector space isomorphism
$V_1\rightarrow V_2$ 
conjugating $W_1$ to $W_2$ (and from which one can obtain an isomorphism with these
properties that preserves the forms, also by an averaging process; see
\cite[\S 14.1]{Kane01}). 
A reflection group is reducible if it has the form $W_1\times W_2$
for non-trivial reflection groups $W_i\subset GL(V)$, and is essential
if only the origin is left fixed by all $g\in W$.

The Shephard-Todd classification (up to this isomorphism) 
of the finite essential irreducible complex
reflection groups then contains three infinite families
and 34 exceptional cases (see for instance, 
\cite[\S 15]{Kane01}). The infinite families are 
the cyclic and symmetric groups, and the groups $G(m,n,p)$ of $n\times n$
monomial matrices whose non-zero entries $\ww_1,\ldots,\ww_n$ are $m$-th roots
of unity with $(\ww_1\ldots\ww_n)^{m/p}=1$.

If $X\subset V$ is a subspace, then
the isotropy group $W_X$ consists of those elements of $W$ that
fix $X$ pointwise. 
Possibly the most significant property of reflection groups for us, 
at least in this paper,
is that $W_X$ is then also a reflection group, generated by reflections
in those hyperplanes containing $X$ \cite{Steinberg60}.

Among the complex groups are the real ones, with the transition from 
a real group $W\subset GL(V_\R)$ to a complex one coming about by passing to
reflections with hyperplanes $H\otimes\C\subset V_\R\otimes\C$.
A finite real reflection group leaves invariant an inner product
$(\,,)$, so that $V$ has the structure of a Euclidean space.

\begin{table}
\begin{tabular}{lll}
\hline
Type and order&Root system $\Phi$&Coxeter symbol and simple system\\\hline
\begin{tabular}{l}
$A_{n-1}\,(n\geq 2)$\\
$n!$\\
\end{tabular}&$\{\xx_i-\xx_j \,\,(1\leq i\not= j\leq n)\}$&$\,\,\,
\begin{picture}(150,30)
\put(-11,0){{$\xx_1-\xx_{2}$}}\put(27,17){{$\xx_{2}-\xx_{3}$}}
\put(68,0){{$\xx_{n-2}-\xx_{n-1}$}}\put(110,17){{$\xx_{n-1}-\xx_{n}$}}
\put(4,9){\circle{8}}\put(8,9){\line(1,0){30}}
\put(42,9){\circle{8}}\put(46,9){\line(1,0){12}}
\put(62,9){$\ldots$}\put(78,9){\line(1,0){12}}
\put(94,9){\circle{8}}\put(98,9){\line(1,0){30}}
\put(132,9){\circle{8}}                 
\end{picture}$\\
\begin{tabular}{l}
$D_n\,(n\geq 4)$\\
$2^{n-1}n!$\\
\end{tabular}&$\{\pm\xx_i\pm\xx_j\,\,(1\leq i<j\leq n)\}$&
\begin{picture}(150,55)
\put(-8,0){{$\xx_1-\xx_{2}$}}
\put(27,17){{$\xx_{2}-\xx_{3}$}}
\put(100,7){{$\xx_{n-2}-\xx_{n-1}$}}\put(105,42){{$\xx_{n-1}-\xx_{n}$}}
\put(105,-27){{$\xx_{n-1}+\xx_{n}$}}
\put(4,9){\circle{8}}\put(8,9){\line(1,0){30}}
\put(42,9){\circle{8}}\put(46,9){\line(1,0){12}}
\put(62,9){$\ldots$}\put(78,9){\line(1,0){12}}
\put(94,9){\circle{8}}                 
\put(96.83,11.83){\line(1,1){21.21}}    
\put(96.83,6.17){\line(1,-1){21.21}}   
\put(120.98,35.98){\circle{8}}          
\put(120.98,-17.98){\circle{8}}         
\end{picture}\\
\begin{tabular}{l}
$B_n\,(n\geq 2)$\\
$2^nn!$\\
\end{tabular}&
\begin{tabular}{l}
$\{\pm\xx_i\,\,(1\leq i\leq n)$,\\
$\pm\xx_i\pm\xx_j\,\,(1\leq i<j\leq n)\}$\\
\end{tabular}
&\begin{picture}(150,50)
\put(-8,0){{$\xx_1-\xx_{2}$}}\put(27,17){{$\xx_{2}-\xx_3$}}
\put(70,0){{$\xx_{n-1}-\xx_{n}$}}\put(128,17){{$\xx_{n}$}}
\put(109,11){$4$}
\put(4,9){\circle{8}}\put(8,9){\line(1,0){30}}
\put(42,9){\circle{8}}\put(46,9){\line(1,0){12}}
\put(62,9){$\ldots$}\put(78,9){\line(1,0){12}}
\put(94,9){\circle{8}}\put(97.5,9){\line(1,0){30}}
\put(132,9){\circle{8}}                 
\end{picture}\\\hline
\end{tabular}\caption{Standard root systems 
$\Phi\subset V$
for the classical Weyl groups
\cite[\S 2.10]{Humphreys90}
where $V$ is a Euclidean space with
orthonormal basis $\{\xx_1,\ldots,\xx_n\}$: 
the type is in the left most column,
with the subscript the dimension of the subspace of $V$ spanned by 
$\Phi$ and the order of the associated Weyl group $W(\Phi)$.
The last column gives the Coxeter symbol, with nodes labelled
by the vectors in a simple system $\Delta\subset\Phi$.}\label{table:roots1}
\end{table}

Traditionally, the finite real groups are studied via
the combinatorics of their root systems: 
an (abstract) 
root system $\Phi$ in a Euclidean space $V$ is a finite
set of non-zero vectors such that, (i). if $\vv\in\Phi$ then
$\lambda\vv\in\Phi$ if and only if $\lambda=\pm 1$, and (ii). if $\uu,\vv\in\Phi$ then
$(\uu)s_\vv\in\Phi$, where $s_\vv$ is the reflection 
in the hyperplane $\vv^\perp$.
The system is essential if the $\R$-span of $\Phi$ is $V$; reducible
if $V=V_1\perp V_2$ and $\Phi=\Phi_1\cup\Phi_2$ for (non-empty) root systems
$\Phi_i\subset V_i$
(in which case we write $\Phi=\Phi_1\perp\Phi_2$),
and crystallographic if
$$
\langle\uu,\vv\rangle:=\frac{2(\uu,\vv)}{(\vv,\vv)}\in\Z,
$$
for all $\uu,\vv\in\Phi$. The associated reflection group is
$W(\Phi)=\langle s_\vv\,(\vv\in\Phi)\rangle$, 
and every finite reflection group arises from
some root system in this way, with the essential, irreducible groups arising
from essential, irreducible systems. The $W(\Phi)$ for $\Phi$ crystallographic
are the {\em Weyl groups\/}.

Root systems $\Phi_i\subset V_i$ are isomorphic if there is an inner product
preserving linear isomorphism $V_1\rightarrow V_2$ sending $\Phi_1$ to $\Phi_2$,
and are {\em stably\/} isomorphic if the isomorphism is between the subspaces
spanned by the $\Phi_i$. In particular, every root system is stably isomorphic to an
essential one. The corresponding groups $W(\Phi_i)$ are stably isomorphic if there
is a vector space isomorphism between the spans of the $\Phi_i$ conjugating one
group to the other.

The irreducible crystallographic root systems have
been classified, up to stable isomorphism:
there are four infinite
families $A,B,C$ and $D$ (the classical systems), 
and five exceptional ones of types $E,F$ and $G$.
The resulting reflection 
groups $W(\Phi)$ provide a list of almost all the finite reflection groups
up to stable isomorphism, with the only omissions being the dihedral groups
and the symmetry groups of the $3$-dimensional dodecahedron/icosahedron and the
$4$-dimensional 120/600-cell.


Table \ref{table:roots1} 
shows the classical crystallographic
$\Phi\subset V$. 
The root systems of types $B$ and $C$ have the same symmetry, 
but different lengths of roots;
nevertheless the associated Weyl groups are identical, and it is these that 
ultimately concern us. We have thus given just the type $B$ system in the
table (type $C$ has roots $\pm 2\xx_i$ rather than the
$\pm\xx_i$). 

\begin{table}
\begin{tabular}{c}
\hline
\begin{pspicture}(0,0)(14.3,2)
\rput(0,-2){
\rput(0,2){
\rput(-.225,-0.2){\rput(3.25,1){\BoxedEPSF{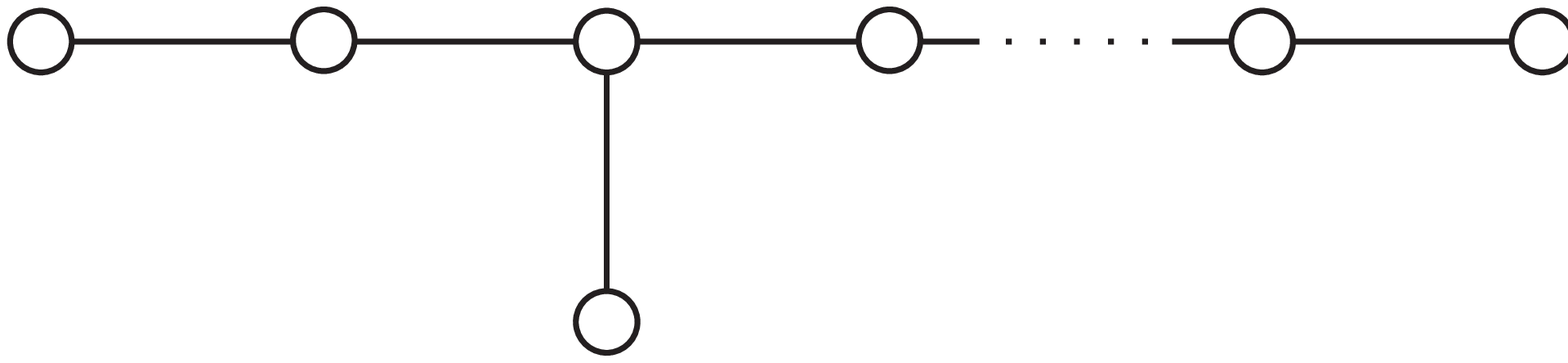 scaled 300}}}
\rput(4,.6){$E_n\,(n=6,7,8)$}
}
\rput(0,-0.2){\rput(12.5,3.3){\BoxedEPSF{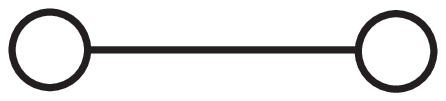 scaled 300}}
\rput(12.5,2.7){$G_2$}}
\rput(12.5,3.3){$6$}
\rput(5.5,-0.2){\rput(3.25,3.3){\BoxedEPSF{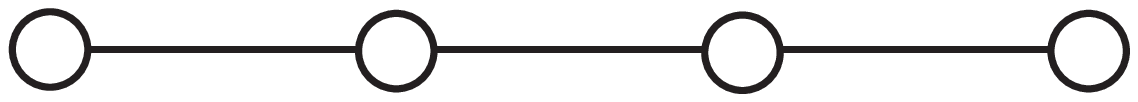 scaled 300}}
\rput(3.2,2.7){$F_4$}\rput(3.25,3.5){$4$}
}}
\end{pspicture}\\
\hline
\end{tabular}
\caption{Symbols for the irreducible exceptional
Weyl groups.}\label{table:roots2}
\end{table}

The last column gives the Coxeter symbol, whose nodes are labelled
by the vectors in a simple system 
$\Delta\subset\Phi$: a basis for the
$\R$-span of $\Phi$ such that each root is a linear combination of $\Delta$ 
with coefficients all of the same sign.
The Weyl group $W(\Phi)$ is then generated by 
the reflections $s_\vv$ for $\vv\in\Delta$
simple. 
The $i$-th and $j$-th nodes of the symbol are connected by an edge
labelled $m_{ij}$, where $\langle\uu,\vv\rangle\langle\vv,\uu\rangle=m_{ij}-2$,
for the simple roots $\uu,\vv$ labelling the nodes, and the rotation
$s_\uu s_\vv$ has order $m_{ij}$ in $W(\Phi)$. It is traditional to omit
labels $m_{ij}=3$, and to remove completely the edges labelled by $m_{ij}=2$.
For convenience in expressing some of the formulae of \S\ref{section:examples},
we adopt the additional 
conventions $A_{-1}=A_0=\varnothing$, $B_0=\varnothing,
B_1=\{\pm\xx_1\}$, and $D_0=D_1=\varnothing$,
$D_n=\{\pm\xx_i\pm\xx_j\,\,(1\leq i<j\leq n)\}$ for $n=2,3$.
In Table \ref{table:roots2} we have given just the Coxeter symbols for 
the exceptional Weyl groups. See \cite[\S 2.10]{Humphreys90} for their root systems.

The Coxeter symbol also gives the \emph{reflectional representation\/} of the Weyl
group: let $S$ be the set of nodes of the symbol and $V$ the real space with basis 
$\{\vv_s\,|\,s\in S\}$ and symmetric bilinear form defined by,
$$
B(\vv_s,\vv_{t})=-\cos\frac{\pi}{m_{st}}.
$$
For $\uu\in V$, define $\sigma_{\uu}:V\rightarrow V$
by $\vv\sigma_{\uu} = \vv-2B(\vv,\uu)\uu$;
then the map $s_{\vv}\mapsto\sigma_{\vv_s}$, where $\vv\in\Delta$ is the label
of $s\in S$, extends to a faithful irreducible representation 
$\sigma:W(\Phi)\rightarrow\gl(V)$.
We will abbreviate $\vv(\sigma(g))$ to $\vv g$.
Any faithful representation of $W(\Phi)$ with the $s_{\vv}$ $(\vv\in\Delta)$
acting as reflections is equivalent to the direct sum of the reflectional 
representation and a trivial representation.

A Weyl group $W(\Phi)$ is of \emph{$(-1)$-type\/} 
if in the reflectional representation there is an element $g\in W(\Phi)$ 
acting on $V$ as the 
antipodal map, ie: $\vv g=-\vv$ for all $\vv\in V$.
They are precisely the groups
with non-trivial center; a reducible Weyl group $W_1\times W_2$ is of $(-1)$-type 
iff each $W_i$ is of $(-1)$-type, 
and the irreducible Weyl groups of $(-1)$-type are the $W(\Phi)$ for
$\Phi=A_1, B_n, D_n (n\text{ even})$, and $E_6$.

It turns out that the classical Weyl groups have alternative descriptions
as certain
permutation groups, and we will use these extensively in this paper.
This is very much in the spirit of the historical development of the
theory of reflection groups, where a number of the classical theorems
were initially proved on a case by case basis, using such descriptions and the 
classification of Coxeter, and while many now have uniform proofs
that intrinsically use the reflection group structure, some still do not.

Firstly then, the map $(i,j)\mapsto s_{\xx_i-\xx_j}$ induces an isomorphism
$\goth{S}_n\rightarrow W(A_{n-1})$, 
and indeed the $W(A_{n-1})$-action on the basis $\{\xx_1,\ldots,\xx_n\}$ is
just permutation of coordinates.

There are two descriptions that prove useful for the 
Weyl group $W(B_n)$. Let 
$I$ be a set and $\QQ(I)$ the collection of all subsets of $I$, which forms an Abelian group
under symmetric difference $X\vartriangle Y:=(X\cup Y)\setminus (X\cap Y)$.
Writing $\prod_I\Z_2$ for the (unrestricted) direct product, we have an 
isomorphism $\prod_I\Z_2\rightarrow\QQ(I)$ given by the map
 $\xx=(x_i)_{i\in I}\mapsto X=\{i\in I\,|\,x_i=1\}$, and this makes it is easy to see that
$\QQ(I)$ is generated by the singletons. The symmetric group $\goth{S}_I$ acts
on $\QQ(I)$ via the obvious $X\mapsto X\sigma$, and thus we may form the 
semi-direct product $\goth{S}_I\ltimes\QQ(I)$, in which every element has a 
unique expression as a pair $\sigma X$, $\sigma\in\goth{S}_I$, $X\subset I$, and
with $\sigma X\tau Y=\sigma\tau (X\tau\vartriangle Y)$.
Write $\goth{S}_n\ltimes\QQ(n)$ if $I=\{1,\ldots,n\}$, in which case the map
$(i,j)\mapsto s_{\xx_i-\xx_j}$, $\{i\}\mapsto s_{\xx_i}$ induces an isomorphism
$\goth{S}_n\ltimes\QQ(n)\rightarrow W(B_n)$.

The second viewpoint is to consider the group 
$\BBB_I$ of {\em signed permutations\/}
of $I$, ie: $\BBB_I=\{\sigma\in\goth{S}_{I\,\cup\,(-I)}\,|\,(-x)\sigma=-(x\sigma)\}$. We then
have an isomorphism $\BBB_n\rightarrow W(B_n)$ induced by
$(i,j)(-i,-j)\mapsto s_{\xx_i-\xx_j}$ and $(i,-i)\mapsto s_{\xx_i}$
(\emph{cf}. Proposition \ref{examples:permutationmonoids:result300}(ii)).

Finally, $\QQ(n)$ has a subgroup $\QQ^+(n)$ consisting of those $X$ with 
$|X|$ even, and the $\goth{S}_n$ action restricting to an action on
$\QQ^+(n)$. The map 
$(i,j)\mapsto s_{\xx_i-\xx_j}$, $\{i,j\}\mapsto s_{\xx_i-\xx_j}s_{\xx_i+\xx_j}$ 
induces an isomorphism $\goth{S}_n\ltimes\QQ^+(n)\rightarrow W(D_n)$. There 
is also a description of $W(D_n)$ in terms of even signed permutations,
but this will be of no use to us.

\section{Systems of subspaces for reflection groups}\label{section:systems}

Partial mirror symmetry describes the phenomenon of restricting the linear 
isomorphisms
of a reflection group to ``local isomorphisms'' between certain subspaces.
In this section we place a modest amount of structure
on these subspaces that still allows for a large number of interesting
examples.

Let $G\subset GL(V)$ be a group. 
A collection $\BB$ of subspaces of $V$ is a {\em system of subspaces for\/}
$G$ if and only if 
\begin{description}
\item [(S1).] $V\in\BB$,
\item [(S2).] $\BB G=\BB$, ie: $Xg\in\BB$ for
any $X\in\BB$ and $g\in G$, and
\item [(S3).] if $X, Y\in\BB$ then $X\cap Y\in\BB$. 
\end{description}

If $\BB_1,\BB_2$ are systems for $G$ then
clearly $\BB_1\cap\BB_2$ is too, and thus for any set $\Omega$ of subspaces 
we write 
$\langle\Omega\rangle_G$ for the intersection of all systems 
for $G$ containing
$\Omega$, and call this the system 
for $G$ {\em generated\/} by $\Omega$.

A system $\BB$ can be partially ordered by inclusion (respectively, reverse inclusion)
and both will turn out to be useful for us. The result is a 
meet (resp. join) semilattice with $\hat{1}$ (resp. $\hat{0}$), 
indeed a lattice if $\BB$ is finite (see \cite[\S 3.1]{Stanley97} for basic 
facts concerning lattices). 
It is an elementary fact in semigroup theory
\cite[Proposition 1.3.2]{Howie95} that a meet semilattice with $\hat{1}$ is
a commutative monoid $E$ of idempotents and vice-versa. For any $e\in E$, let
$Ee=\{x\in E\,|\,x\leq e\}$. The Munn semigroup 
\cite[\S 5.4]{Howie95} $\TT_E$ of $E$
is then defined to be the set of all isomorphisms $Ee\rightarrow Ef$
where $e,f$ range over all elements of $E$ with. The following is then easily proved:

\begin{proposition}\label{systems:result50}
$\BB$ is a system in $V=\F\text{-span}\{\BB\}$ for $G\subset GL(V)$ if and only if
$E=(\BB,\cap)$ is a commutative monoid of idempotents and the mapping
$g\mapsto \theta_g$ where $X\theta_g = Xg$ for $X\in \BB$ and $g\in G$ 
 is a (monoid)
homomorphism $G\rightarrow\TT_E$ to the Munn semigroup of $E$.
\end{proposition}


Recall 
that in a poset $(\PP,\leq)$, if
$x<y$ and there is no
$z$ with $x<z<y$ then we say that $y$ covers $x$, and write
$x<_c y$. $\PP$
is graded of rank $n$
if every chain $x_1<_c\cdots<_c x_n$,
maximal under inclusion of such chains,
has the same length $n$.
There is then a unique rank function $\rk:\PP\rightarrow\{0,1,\ldots,n\}$ with $\rk(x)=0$ 
if and only if $x$ is minimal, and 
$\rk(y)=\rk(x)+1$ whenever $x<_cy$.
The rank $1$ (resp. rank $n-1$) elements of $\PP$ are
the atoms (resp. coatoms) and $\PP$ is atomic (resp. coatomic) if every
element is a join of atoms (resp. meet of coatoms).
A Boolean lattice on a finite set $X$
is a lattice isomorphic to the lattice of all subsets of $X$ under 
inclusion.

In particular, if we order a 
system $\langle\Omega\rangle_G$ of subspaces 
for $G\subset GL(V)$ by inclusion (resp. reverse inclusion),
then every element is a meet (resp. join) of
the $Xg$ for $X\in\Omega$ and $g\in G$; if $V$ is finite dimensional
and all the $X\in\Omega$ have the 
same dimension, then we have 
a coatomic poset with coatoms (resp. atomic poset with atoms) the $Xg$.

A \emph{hyperplane arrangement\/} $\AA$ is a finite collection of hyperplanes
in $V$. General references are \cite{Orlik92,Zaslavsky75},
where the hyperplanes are allowed to be affine, but we will restrict ourselves to 
arrangements where the hyperplanes are linear (hence subspaces of $V$).
An important combinatorial invariant for $\AA$ is the
intersection lattice $L(\AA)$--the set of all possible intersections of elements
of $\AA$, ordered by reverse inclusion, and with the null intersection taken to be
the ambient space $V$. What results is a graded atomic lattice 
of rank $\codim \bigcap_{X\in\AA} X$ \cite[\S 2.1]{Orlik92},
with $\hat{0}$ the space $V$, atoms the hyperplanes in $\AA$
and $\rk X=\codim X$.

If $G\subset GL(V)$ is finite and $\AA\subset V$ a hyperplane arrangement, then
$\AA G$ is also a hyperplane arrangement, for which the following is
then obvious,

\begin{lemma}\label{systems:result100}
The system $\langle\AA\rangle_G$ for $G$ generated by $\AA$ is the intersection 
lattice $L(\AA G)$, and the $G$-action on $\langle\AA\rangle_G$ is
rank preserving.
\end{lemma}

In general $L(\AA)\subset L(\AA G)$, but we will often have 
$\AA G=\AA$, hence equality of the intersection lattices.


\subsection{Boolean systems}
\label{subsection:permutationsystem}

Specializing now to reflection groups, 
a simple but nevertheless interesting example of a system arises if
$V$ is a Euclidean space 
with orthonormal basis $\{\xx_1,\ldots,\xx_n\}$ and 
$W=W(\Phi)$ a Weyl group
as in Tables \ref{table:roots1}-\ref{table:roots2}. 
The {\em Boolean\/} (or {\em orthogonal\/}) hyperplane arrangement 
\cite[\S 1.2]{Orlik92}
$\AA=\{\xx_1^\perp,\ldots,\xx_n^\perp\}$ consists of the coordinate hyperplanes,
and we call the system $\langle\AA\rangle_W$ for $W$ generated by $\AA$ a
{\em Boolean system}. The name stems from the fact that $L(\AA)$ is a Boolean
lattice, although it should be noted that the system $\langle\AA\rangle_W$ 
itself will only be Boolean when we have $\AA W=\AA$. 

Consider a Weyl group  $W=W(\Phi)$ with $\Phi$ a classical root system
as in Table \ref{table:roots1}. Then $\AA W=\AA$, and hence
the Boolean system $\langle\AA\rangle_W=L(\AA)$ is a Boolean lattice with
the map 
$\xx_{i_1}^\perp\cap\cdots\cap \xx_{i_k}^\perp\mapsto\{i_1,\ldots,i_k\}$
being a lattice isomorphism from $L(\AA)$ to the lattice of subsets of 
$I=\{1,2,\ldots,n\}$.

%

The rank $k$ elements of $L(\AA)$ are the intersections 
$\xx_{i_1}^\perp\cap\cdots\cap\xx_{i_k}^\perp$ of $k$ distinct hyperplanes, and as
the symmetric group $\goth{S}_X$ for $X=\{\xx_1,\ldots,\xx_n\}$ is a subgroup of $W(\Phi)$
for classical $\Phi$, 
the action of $W(\Phi)$ on the rank $k$ elements is transitive.

If $\Phi$ is a root system for one
of the exceptional groups in Table \ref{table:roots2},
then $\AA\subset\AA W$, but the system for $W(\Phi)$
will have more elements than the intersection lattice $L(\AA)$. 
If for instance $\Phi$ is the $F_4$ root system of \cite[\S 2.10]{Humphreys90}
and $W=W(\Phi)$, then the system
$\langle\AA\rangle_W$ has atoms the hyperplanes $\xx_i^\perp$ and 
$\frac{1}{2}(\pm\xx_1\pm\xx_2\pm\xx_3\pm\xx_4)^\perp$, ie: the reflecting hyperplanes
of $W(F_4)$ corresponding to the short roots, and as such is a
subsystem of the intersection lattice of the type $F_4$ reflection arrangement
of \S \ref{subsection:systemintersectionposet}. If $\Phi=E_6,E_7$ or $E_8$,
then a description of the Boolean system is possible, but messier.

\subsection{Intersection lattices of reflection arrangements} 
\label{subsection:systemintersectionposet}

A more natural example of a system of subspaces for a reflection group
$W$ is given by the intersection lattice $L(\AA)$ of the 
reflecting hyperplanes $\AA$ of $W$. If $X\in\AA$ and $s_X^{}\in W$ is the
reflection in $X$, then for $g\in W$ we have $s_{Xg}^{}=g^{-1}s_X^{}g$, and 
so $Xg\in\AA$. Thus $\AA W=\AA$, and we have,

\begin{lemma}\label{systems:result200}
If $W\subset GL(V)$ is a reflection group and $\AA$ the hyperplane arrangement
consisting of the reflecting hyperplanes of $W$, then 
$\langle\AA\rangle_W=L(\AA)$.
\end{lemma} 

We will call such an $L(\AA)$ a {\em (reflection) arrangement system\/}, and
for the remainder of this section we focus on these systems 
(ordered by reverse inclusion) when $W$ is
a Weyl group as in Tables \ref{table:roots1}-\ref{table:roots2},
summarizing the necessary results of
\cite[\S 6.4]{Orlik92}.
Recall that a partition of $I=\{1,2,\ldots,n\}$ is a collection
$\Lambda=\{\Lambda_1,\ldots,\Lambda_p\}$ of nonempty pairwise disjoint subsets
$\Lambda_i\subset I$ whose union is $I$. If $\lambda_i=|\Lambda_i|$ then 
$\lambda=\|\Lambda\|=(\lambda_1,\ldots,\lambda_p)$ is a partition of $n$, ie: 
the integers $\lambda_i\geq 1$ with $\sum\lambda_i=n$, and we 
order the $\Lambda_i$ so that
$\lambda_1\geq\cdots\geq\lambda_p\geq 1$. 
Order the set $\Pi(n)$ of partitions of $I$ by
refinement, ie: $\Lambda\leq\Lambda'$
if and only if for every $\Lambda_i$ there is a $\Lambda'_j$ with $\Lambda_i\subset\Lambda'_j$.
The result is an atomic graded lattice with 
$\rk\Lambda=\sum(\lambda_i-1)$ and atoms
the $\Lambda$ with $\lambda_1=2$ and $\lambda_i=1$ for $i>1$. The following
is \cite[Proposition 2.9]{Orlik92}:

\begin{proposition}
\label{systems:result300}
Let $\AA$ be the hyperplane arrangement consisting of the reflecting hyperplanes 
of the Weyl group $W(A_{n-1})$.
Then the map that sends the atomic partition with $\Lambda_1=\{i,j\}$
to the hyperplane $(\xx_i-\xx_j)^\perp$ extends to 
a lattice isomorphism $\Pi(n)\rightarrow L(\AA)$.
\end{proposition}

Indeed, writing $X(\Lambda)\in L(\AA)$ for the image of $\Lambda$, we have
\begin{equation}\label{subspaces:typeA}
X(\Lambda)=\bigcap_{\lambda_k>1}\kern4pt\bigcap_{i,j\in\Lambda_k}(\xx_i-\xx_j)^\perp.
\end{equation}
For a partition $\Lambda$, let $b_i>0$ be the number of
$\lambda_j$ equal to $i$, and 
$$
b_\lambda=b_1!b_2!\ldots (1!)^{b_1}(2!)^{b_2}\ldots
$$ 
If $\sigma\mapsto g(\sigma)$ is the isomorphism
$\goth{S}_n\rightarrow W(A_{n-1})$ of \S \ref{section:reflectiongroups},
then the action of $W(A_{n-1})$ on $L(\AA)$ is given by 
$X(\Lambda)g(\sigma)=X(\Lambda\sigma)$, where 
$\Lambda\sigma=\{\Lambda_1\sigma,\ldots,\Lambda_p\sigma\}$.
The following is \cite[Proposition 6.72]{Orlik92}:

\begin{proposition}
\label{systems:result400}
In the action of the Weyl group $W(A_{n-1})$ on %
$L(\AA)$,
two subspaces $X(\Lambda)$ and $X(\Lambda')$ lie in the same orbit if
and only if 
$\|\Lambda\|=\|\Lambda'\|$. The cardinality of the orbit of
the subspace $X(\Lambda)$ is $n!/b_\lambda$.
\end{proposition}

Turning now to the Weyl group $W(B_n)$, 
let $\TT(I)$ be the set of triples $(\Delta,\Gamma,\Lambda)$
where $\Delta\subset I$, $\Gamma\subset J:=I\setminus\Delta$
and $\Lambda=\{\Lambda_1,\ldots,\Lambda_p\}$ 
is a partition of $J$.
There is then \cite[Proposition 6.74]{Orlik92} a surjective mapping
$\TT(I)\ni(\Delta,\Gamma,\Lambda)\mapsto 
X(\Delta,\Gamma,\Lambda)\in L(\AA)$,
with 
\begin{equation}\label{subspaces:typeB}
X(\Delta,\Gamma,\Lambda)
=\bigcap_{\lambda_k>1}\kern3pt\bigcap_{i,j\in\Lambda_k}(\xx_i+\ve_j\xx_j)^\perp
\cap\bigcap_{i\in\Delta}\xx_i^\perp
\,\,
\text{ and }
\,\,
\rk X(\Delta,\Gamma,\Lambda)=|\Delta|+\sum(\lambda_i-1),
\end{equation}
where $\ve_j=1$ if $j\in\Gamma$ or $\ve_j=-1$ if $j\not\in\Gamma$.
Moreover,
$X(\Delta,\Gamma,\Lambda)=X(\Delta',\Gamma',\Lambda')$ 
if and only if $\Delta=\Delta'$,
$\Lambda=\Lambda'$ and for each 
$1\leq i\leq p$, $\Gamma_i'=\Gamma_i$ or $\Lambda_i\setminus\Gamma_i$,
where $\Gamma_i=\Gamma\cap\Lambda_i$ and $\Gamma'_i$ is defined similarly.

If $\sigma T\mapsto g(\sigma,T)$ is the isomorphism 
$\goth{S}_n\ltimes\QQ(n)\rightarrow W(B_n)$ of 
\S\ref{section:reflectiongroups}, then the action of $W(B_n)$
on $L(\AA)$ is given by
$$
X(\Delta,\Gamma,\Lambda)g(\sigma,T)=X(\Delta\sigma,(T_J\vartriangle\Gamma)\sigma,\Lambda\sigma),
$$
where $T_J=T\cap J$.

\begin{proposition}
\label{systems:result500}
In the action of the Weyl group $W(B_n)$ on 
$L(\AA)$, 
two subspaces $X(\Delta,\Gamma,\Lambda)$ and $X(\Delta',\Gamma',\Lambda')$ 
lie in the 
same orbit if and only if $|\Delta|=|\Delta'|$ and 
$\|\Lambda\|=\|\Lambda'\|$. The cardinality of the orbit of
the subspace $X(\Delta,\Gamma,\Lambda)$ is 
$$
2^{j-p}\left(\begin{array}{c}n\\j\end{array}\right)\frac{j!}{b_\lambda},
$$
where $j=|J|$ and $\Lambda=\{\Lambda_1,\ldots,\Lambda_p\}$.
\end{proposition}

(See \cite[Proposition 6.75]{Orlik92}.
What Orlik and Terao actually describe 
is the corresponding result for the full monomial group $G(r,1,n)$, 
where we have contented ourselves with $G(2,1,n)\cong W(B_n)$.)
For the Weyl group $W(D_n)$ and its reflecting hyperplanes $\AA$, 
let $\SSS(I)$ be the subset of $\TT(I)$
consisting of those triples $(\Delta,\Gamma,\Lambda)$ with $|\Delta|\not= 1$.
Then by \cite[Proposition 6.78]{Orlik92} there is a surjective mapping 
$\SSS(I)\rightarrow L(\AA)$,
where 
\begin{equation}\label{subspaces:typeD}
X(\Delta,\Gamma,\Lambda)
=\left\{\begin{array}{l}
{\displaystyle \bigcap_{\lambda_k>1}\kern4pt\bigcap_{i,j\in\Lambda_k}
(\xx_i+\ve_j\xx_j)^\perp
\text{ if }\Delta=\varnothing},\\
{\displaystyle \bigcap_{\lambda_k>1}\kern4pt\bigcap_{i,j\in\Lambda_k}
(\xx_i+\ve_j\xx_j)^\perp
\cap\bigcap_{i,j\in\Delta}(\xx_i+\xx_j)^\perp\cap(\xx_i-\xx_j)^\perp,
\text{ if }|\Delta|\geq 2},
\end{array}\right.
\end{equation}
and $X(\Delta,\Gamma,\Lambda)=X(\Delta',\Gamma',\Lambda')$ if and only if 
$\Delta=\Delta'$,
$\Lambda=\Lambda'$ and for each 
$1\leq i\leq p$, $\Gamma_i'=\Gamma_i$ or $\Lambda_i\setminus\Gamma_i$,
where $\Gamma_i=\Gamma\cap\Lambda_i$ and $\Gamma'_i$ is defined similarly.
\cite[Proposition 6.79]{Orlik92} then gives,

\begin{proposition}
\label{systems:result600}
If $X(\Delta,\Gamma,\Lambda)$ and $X(\Delta',\Gamma',\Lambda')$ lie in the 
same orbit of the action of $W(D_n)$ on 
$L(\AA)$, then
$|\Delta|=|\Delta'|$ and $\|\Lambda\|=\|\Lambda'\|$. Conversely, suppose 
that $|\Delta|=|\Delta'|$ and $\|\Lambda\|=\|\Lambda'\|$.
\begin{enumerate}
\item If $|\Delta|\geq 2$ then $X(\Delta,\Gamma,\Lambda)$ and 
$X(\Delta',\Gamma',\Lambda')$ lie in the same orbit, which has cardinality 
as in Proposition \ref{systems:result500}.
\item If $\Delta=\varnothing$, then the $W(B_n)$ orbit 
determined by $\|\Lambda\|=(\lambda_1,\ldots,\lambda_p)$ 
forms a single $W(D_n)$ orbit, except when each $\lambda_i$ is even,
in which case it decomposes
into two $W(D_n)$ orbits of size
$$
\frac{2^{n-p-1}\,n!}{b_\lambda}.
$$
\end{enumerate}
\end{proposition}

In part 2 of Proposition \ref{systems:result600}, and when all the $\lambda_i$ are 
even, one of the $W(D_n)$ orbits consists of the 
$X(\varnothing,\Gamma,\Lambda)$ with $|\Gamma|$ even,
and the other with the $|\Gamma|$ odd
(again, Orlik and Terao deal with the monomial group $G(r,r,n)$, while we
consider only $G(2,2,n)\cong W(D_n)$, with the
decomposition of the second part of Proposition \ref{systems:result600} being 
into $d$ $W(D_n)$-orbits, for $d$ the greatest common divisor of 
$\{r,\lambda_1,\ldots,\lambda_p\}$).

\begin{table}
\begin{tabular}{cl}
\hline
{\small{\tt g2}}&
\begin{tabular}{l}
{\small{\tt 1a0:3a1.3a1:1g2}}\\
\end{tabular}\vrule width 0 mm height 0 mm depth 2 mm\\
{\small{\tt f4}}&
\begin{tabular}{l}
{\small{\tt 1a0:12a1.12a1:72a12.16a2.16a2.18b2:12b3.12b3.48a1a2.48a1a2:1f4}}\\
\end{tabular}\vrule width 0 mm height 0 mm depth 2 mm\\
{\small{\tt e6}}&
\begin{tabular}{l}
{\small{\tt 1a0:36a1:270a12.120a2:540a13.720a1a2.270a3:1080a12a2.120a22}}\\
{\small{\tt 540a1a3.216a4.45d4:360a1a22.216a1a4.36a5.27d5:1e6}}\\
\end{tabular}\vrule width 0 mm height 0 mm depth 4 mm\\
{\small{\tt e7}}&
\begin{tabular}{l}
{\small{\tt 1a0:63a1:945a12.336a2:315a13.3780a13.5040a1a2.1260a3:3780a14}}\\
{\small{\tt 15120a12a2.3360a22.1260a1a3.7560a1a3.2016a4.315d4:5040a13a2}}\\
{\small{\tt 10080a1a22.7560a12a3.5040a2a3.6048a1a4.336a5.1008a5.945a1d4}}\\
{\small{\tt 378d5:5040a1a2a3.2016a2a4.1008a1a5.288a6.378a1d5.63d6.28e6:1e7}}\\
\end{tabular}\vrule width 0 mm height 0 mm depth 8 mm\\
{\small{\tt e8}}&
\begin{tabular}{l}
{\small{\tt a1a0:120a1:3780a12.1120a2:37800a13.40320a1a2.7560a3:113400a14}}\\
{\small{\tt 302400a12a2.67200a22.151200a1a3.24192a4.3150d4:604800a13a2}}\\
{\small{\tt 403200a1a22.453600a12a3.302400a2a3.241920a1a4.40320a5.37800a1d4}}\\
{\small{\tt 7560d5:604800a12a22.604800a1a2a3.362880a12a4.151200a32.241920a2a4}}\\
{\small{\tt 120960a1a5.34560a6.50400a2d4.45360a1d5.3780d6.1120e6:241920a1a2a4}}\\
{\small{\tt 120960a3a4.34560a1a6.8640a7.30240a2d5.1080d7.3360a1e6.120e7:1e8}}\\
\end{tabular}\\
\hline
\end{tabular}
\caption{Orbit data for the exceptional arrangement systems
\cite[Appendix C]{Orlik92}:
each orbit is encoded in a string consisting of the 
number of subspaces in the orbit followed by their common stabilizer
written in the form {\tt xnmypq...}, to indicate the product of Weyl groups
$X_n^m\times Y_p^q\ldots$
Different orbits
of subspaces of the same rank are separated by a period and orbits of different
ranks by a colon.}\label{table:orbitdata}
\end{table}

If $W$ is an exceptional Weyl group then a convenient description
of $L(\AA)$ is harder, but an enumeration of the orbits
of the $W$-action on $L(\AA)$ suffices for our purposes. We summarize some of the
results of \cite{Orlik83,Orlik82} (see \cite[Appendix C]{Orlik92}) in Table
\ref{table:orbitdata}. For
example, the orbit data for the Weyl group $W(E_6)$, which starts as,
\begin{center}
{\tt 1a0:36a1:270a12.120a2:540a13.720a1a2.270a3}
\end{center}
indicates a single rank $0$ orbit with stabilizer the Weyl group $A_0\cong 1$
(corresponding to the ambient space $V$), a single rank $1$
orbit of size $36$ with stabilizer $A_1\cong\Z_2$ (corresponding to the 
reflecting hyperplanes, or
the $72$ roots in the $E_6$ root system arranged in $36$ $\pm$ pairs),
two orbits of
rank $2$ subspaces of sizes $270$ and $120$ with stabilizers $A_1\times A_1$ and
$A_2$ respectively, and so on.
There are distinct rank one orbits with isomorphic stabilizers in types $G_2$ and
$F_4$, corresponding to the two conjugacy classes of generating reflections
(this phenomenon not arising in type $E$ where all the generating reflections
are conjugate).

We have stuck to the Weyl groups, as promised in \S \ref{section:reflectiongroups}, 
but the data
in Table \ref{table:orbitdata} could just as easily be read off 
\cite[Appendix C]{Orlik92} for all 34 exceptional 
finite complex reflection groups.

\section{Inverse Monoids and Reflection Monoids}\label{section:reflectionmonoids}

We are now ready for reflection monoids and some of
their elementary properties, but first we recall some of the basic concepts
of inverse monoids. For more on the general theory of inverse monoids
see \cite[Chapter 5]{Howie95} and \cite{Lawson98}.

An \emph{inverse monoid} is a monoid $M$ such that for all $a\in M$
there is a unique $b\in M$ such that $aba = a$ and $bab = b$. The
element $b$ is the \emph{inverse} of $a$ and is denoted by $a^{-1}$.
It is worth noting that $(a^{-1})^{-1} =a$ and $(ab)^{-1} =
b^{-1}a^{-1}$ for all $a,b \in M$.
The set of idempotents $E(M)$ of $M$ forms a commutative submonoid,
referred to as the \emph{semilattice of idempotents} of $M$. We
denote the group of units of $M$ by $G(M)$.
An \emph{inverse submonoid} of an inverse monoid $M$ is simply a
submonoid $N$
closed under taking inverses; it is \emph{full} if $E(N) = E(M)$.

The archetypal example of an inverse monoid is the symmetric inverse
monoid defined as follows. 
For a non-empty set $X$, a \emph{partial
permutation} is a bijection $\sigma:Y\to Z$ for some subsets $Y,Z$ of $X$.
We allow $Y$ and $Z$ to be empty so that the empty function is regarded
as a partial permutation. The set of all partial permutations of $X$ is made
into a monoid by using the usual rule for composition of partial
functions; it is called the \emph{symmetric
inverse monoid} on $X$ and denoted by $\II_X$ (if $X =
\{1,2,\dots,n\}$, we write $\II_n$ for $\II_X$). That it is an inverse
monoid follows from the fact that if $\sigma$ is a partial permutation of
$X$, then so is its inverse  (as a function) $\sigma^{-1}$, and this is the
inverse of $\sigma$ in $\II_X$ in the sense above. 
Clearly, the group of
units of $\II_X$ is the symmetric group $\goth{S}_X$, and $E(\II_X)$ consists
of the \emph{partial identities} $\ve_Y^{}$ for all subsets $Y$ of $X$
where $\ve_Y^{}$ is the identity map on the subset $Y$.
It is clear that, for $Y,Z \subset X$, we have $\ve_Y^{}\ve_Z^{} = \ve_{Y\cap
Z}^{}$ and hence that $E(\II_X)$ is isomorphic to the Boolean algebra of
all subsets of $X$.  

Just as $\goth{S}_n$ is isomorphic
to the group of permutation matrices, so $\II_n$ is isomorphic to the 
monoid of partial permutation matrices, or {\em rook monoid}:
the $n\times n$ matrices having $0,1$ entries with at most one non-zero entry
in each row and column (and so called as each element represents an $n\times n$
chessboard with the $0$ squares empty and the $1$ squares containing rooks,
with the rooks mutually non-attacking).

We observe that if $M$ is an inverse submonoid of
$\II_X$, then
$$E(M) = M \cap  E(\II_X) = \{ \ve_Y^{}\mid Y = \dom\, \sigma \text{ for some } 
\sigma \in M\}.$$
Equally, $E(M) = \{ \ve_Y^{}\mid Y = \im\, \sigma \text{ for some } \sigma \in
M\}$ since $\im \sigma = \dom\, \sigma^{-1}$ for all $\sigma \in M$.
Putting
$$\BB = \{\, \dom\, \sigma \mid \sigma \in M\},$$
we see that $\BB$ is a meet semilattice isomorphic to $E(M)$. Moreover,
$X\in \BB$ since $M$ is a submonoid, and finally, if $Y\in \BB$ and 
$g \in G(M)$, then $Yg = \im (\ve_{Y}^{}g) \in
\BB$. Thus $\BB$ satisfies analogues of 
(S1)-(S3) in 
\S\ref{section:systems} for a system of subspaces for a subgroup
of $ GL(V)$, so we say that it is a \emph{system of subsets} for the
group $G(M)$. 

Every inverse monoid $M$
has a faithful representation (called the Vagner-Preston representation)
$\rho_M:M\to \II_M$ by partial permutations given by partial right
multiplication \cite{Howie95,Lawson98}, and the significance of
the symmetric inverse monoid is due partly to this fact.

Another example of an inverse monoid that we will encounter in
\S\ref{subsection:permutationmonoids}
is the {\em monoid of partial signed permutations\/} of a non-empty
set $X$. Let $-X = \{-x\mid x\in X\}$ be disjoint from $X$ such that $x
\mapsto -x$ is a bijection, and
define
$$
\JJ_X:=\{\sigma\in\II_{X\cup-X}\,|\,(-x)\sigma=-(x\sigma)\text{ and }x\in\dom\,\sigma
\Leftrightarrow -x\in\dom\,\sigma\},
$$
where we write $\JJ_n$ when $X=\{1,2,\ldots,n\}$ and in this case $-x$
has its usual meaning. The group of units of $\JJ_X$ is the group 
$\BBB_X$ of partial signed permutations of $X$.


We shall be particularly interested in factorizable inverse monoids,
where an inverse monoid $M$ is \emph{factorizable} if $M = E(M)G(M)\ (=
G(M)E(M))$. See \cite{Lawson98} for more details regarding  
factorizable inverse monoids. For $\sigma \in M$ 
where $M$ is an inverse submonoid of $\II_X$, we have
$\sigma \in E(M)G(M)$ if and only if $\sigma$ is a restriction of a unit of
$M$, so that factorizable inverse submonoids of $\II_X$ are those in
which every element is a restriction of some unit of $M$. For example,
$\II_n$ is factorizable, since any partial permutation of
$\{1,\dots,n\}$ can be extended (not necessarily uniquely) to an element
of $\goth{S}_n$. 
However, if $X$ is infinite, then $\II_X$ is not
factorizable since, for example, an injective map from $X$ to itself
(with domain $X$) which is not surjective cannot be a restriction of a
permutation of $X$. Similarly, $\JJ_n$ is factorizable, but $\JJ_X$ is
not when $X$ is infinite.

Let $\BB$ be a system of subsets for a subgroup $G$ of $\goth{S}_X$ and
define 
$$F = M(G,\BB) = \{ g_Y^{} \mid g\in G,\, Y\in \BB\}$$
where $g_Y^{}$ is the restriction of $g$ to the subset $Y$. Note that $F
\subset \II_X$ and that if $g_Y^{},h_Z^{} \in F$, then $(g_Y^{})^{-1} =
(g^{-1})_{Yg}^{} \in F$ and $g_Y^{}h_Z^{} = (gh)_T^{}$ 
with $T = Y \cap Zg^{-1}$, so
that $F$ is an inverse submonoid of $\II_X$. Clearly, $G$ is the group
of units of $F$, and $E(F) = \{ \ve_Y^{} \mid Y\in \BB\}$. Moreover, every
element of $F$ is a restriction of a unit, so $F$ is factorizable.

Now let $M$ be an inverse submonoid of $\II_X$ 
and $G$ be its group of units.
Let $\BB$ be the system of subsets for $G$ described above, that is, 
$$ \BB = \{\, \dom\, \sigma \mid \sigma \in M\}.$$ 
Put $F_M = M(G,\BB)$ and note that $F_M$ is a factorizable inverse
submonoid of $M$, and, in fact, it is the largest such submonoid
(\emph{cf.} \cite[Proposition~2.2.1]{Lawson98}). 

Thus if $M$ is actually factorizable, then $M = F_M$, and since every
inverse monoid can be embedded in some $\II_X$, we have a description of
all factorizable inverse monoids. As an illustration, we note that
$\II_n$ can be realised as $M(\goth{S}_n,\BB)$ where $\BB$ is the power set
of $\{1,\dots,n\}$.

Another class of inverse monoids of interest to us are the
fundamental inverse monoids. On any inverse monoid $M$, define the
relation $\mu$ by the rule:
$$ a\,\mu\, b \text{ if and only if }  a^{-1}ea = b^{-1}eb \text{ for all
} e\in E(M).$$
It is easy to see that $\mu$ is a congruence on $M$; it is
idempotent-separating in the sense that distinct idempotents in $M$ are
not related by $\mu$, and, in fact, it is the greatest
idempotent-separating congruence on $M$. We say that $M$ is
\textit{fundamental} if $\mu$ is the equality relation, and mention that
for any $M$, the monoid $M/\mu$ is fundamental. The Munn semigroup
$\TT_E$ of a
semilattice $E$ that we introduced in \S 2 plays a crucial role in
describing fundamental inverse monoids. First, we note that $\TT_E$ is
an inverse submonoid of $\II_E$ whose semilattice of idempotents is
isomorphic to $E$ (see \cite[Theorem 5.4.4]{Howie95} or \cite[Theorem 5.2.7]{Lawson98}). 

Given any inverse monoid $M$ and $a\in M$, define an element $\delta_a
\in \TT_{E(M)}$ as follows. The domain of $\delta_a$ is $Eaa^{-1}$ and
$x\delta_a = a^{-1}xa$ for $x\in Eaa^{-1}$. Note that $\im \delta_a = Ea^{-1}a$.
The main results are the following, for which one should
consult
\cite[Theorem 5.4.4]{Howie95}, \cite[Theorem 5.2.8]{Lawson98} and
\cite[Theorem 5.4.5]{Howie95}, \cite[Theorem 5.2.9]{Lawson98}.

\begin{proposition}
\label{inverse_monoids:result50}
If $M$ is an inverse monoid, then the mapping $\delta:M \to \TT_{E(M)}$
given by $a\delta = \delta_a$ is a homomorphism onto a full inverse
submonoid of $\TT_{E(M)}$ such that $a\delta = b\delta$ if and only if
$a\,\mu\,b$.
\end{proposition}

\begin{proposition} 
\label{inverse_monoids:result100}
An inverse monoid $M$ is fundamental if and only if $M$ is isomorphic to
a full inverse submonoid of $\TT_{E(M)}$.
\end{proposition}

The homomorphism $\delta:M \to \TT_{E(M)}$ of
Proposition~\ref{inverse_monoids:result50} is called the
\textit{fundamental} or \textit{Munn
representation} of $M$. Note that $M$ is fundamental if and only if
$\delta$ is one-one.
 
It is well known that $\II_X$ is fundamental for any set $X$ (see, for
example, \cite[Chapter 5, Exercise~22]{Howie95}). In contrast, for any
nonempty set $X$, it is easy to see that $\JJ_X$ is
\textit{not} fundamental: a simple calculation shows that the identity
of $\JJ_X$ and the transposition $(x,-x)$ are $\mu$-related. In the next
section we see that $\JJ_n$ is a reflection monoid, so there are
non-fundamental reflection monoids. 

We now describe fundamental factorizable inverse monoids  in terms of 
semilattices and their automorphism
groups, a point of view that will prove useful in 
\S \ref{subsection:rennermonoids}. We
remark that the principal ideals of a semilattice $E$ regarded as a monoid are
precisely the principal order ideals of $E$ regarded as a partially
ordered set. It will be convenient to write $\ve_x$ for the partial
identity with domain $Ex$. 

\begin{proposition}\label{reflection_monoids:result50}
If $E$ is a semilattice with greatest element
$\hat{1}$ and $G$ is a subgroup of the automorphism group $\aut(E)$,
then the collection 
$$\BB = \{ Ex \mid x\in E\}$$
 of all principal
ideals of $E$ forms a system of subsets (of $E$) for $G$,
 and
the resulting $M(G,\BB)$ is 
the submonoid of $\TT_E$ generated by $G$ and $E$.

Conversely, any fundamental factorizable inverse monoid $M$ is isomorphic
to a submonoid of $\TT_{E(M)}$ generated by a group $G$ of automorphisms
of $E(M)$ and $E(M)$.
\end{proposition} 
\begin{proof}
Given $E$ and $G$ we observe that 
 $\BB$ does  form a system of subsets (in $E$) for $G$ since $E=
E\hat{1}$, $Ex \cap Ey = Exy$ and the image under $g\in G$ of $Ex$ is
$E(xg)$. 
We can thus define 
the factorizable inverse monoid $M(G,\BB)\subset\II_E$ as above. As $G$
is a subgroup of $\aut(E)$, it is a subgroup of the group of units of
$\TT_E$, and hence if $\ve_xg \in M(G,\BB)$ with $g\in G$, then $\ve_xg
\in \TT_E$. Thus $M(G,\BB)\subset \TT_E$; in fact, it is clearly a full
inverse submonoid of $\TT_E$ and so it is fundamental. Identifying $E(\TT_E)$
with $E$, it is also
clear that $M(G,\BB)$ is generated as a submonoid by $G$ and $E$.

For the converse, let $F$ be a fundamental factorizable inverse monoid
and write $E$ for $E(F)$.
Then $F$ is isomorphic to to a full submonoid of $\TT_E$ which we
identify with $F$. The group $G = G(F)$ of units of $F$ is a 
subgroup of the group
of units of $\TT_E$, that is, of $\aut(E)$. As above $\BB = \{\dom\,\sigma
\mid \sigma \in F\}$ is a system of subsets (of $E$) for $G$ and since $F$
is factorizable, $F = M(G,\BB)$. Thus $F$ is generated by $G$ and $E$
(identifying $E$ with $E(\TT_E)$).
\qed
\end{proof}

If the semilattice $E$ has a least element $\hat{0}$ (in particular, if
$E$ is a lattice), then the principal ideal $Ex$ of $E$ is just the
interval $[\hat{0},x] =\{z\in E\,|\,\hat{0}\leq z\leq x\}$, so that the
system described above is the collection of intervals $\BB=\{[\hat{0},x]\,|\,x\in E\}$.


We now turn to reflection monoids. 
Throughout the rest of this section, $V$ is a vector space over a field $\F$.
A {\em partial linear isomorphism\/} of $V$ is a vector space isomorphism 
$\aa:X\rightarrow Y$ between vector subspaces $X,Y$ of $V$. Thus the set
$ ML(V)$ of all partial isomorphisms of $V$ is a subset of $\II_V$. In
fact, it is an inverse submonoid of $\II_V$ since the composition of two
partial isomorphisms is easily seen to be a partial isomorphism, and the
inverse of an isomorphism is again an isomorphism. The group of units of
$ ML(V)$ is $ GL(V)$, and the semilattice of idempotents consists of all
the partial identities on subspaces of $V$. 

If $V$ has finite
dimension and $X$ is a subspace, then by extending a basis of $X$, 
any partial isomorphism with domain $X$ can be extended to a (not
necessarily unique) full isomorphism of $V$. Thus every element of
$ ML(V)$ is a restriction of a unit, so that $ ML(V)$ is factorizable.
Of course, this is not the case if $V$ has infinite dimension. We record
these observations in the next result.

\begin{lemma}\label{reflection_monoids:result100}
The set $ ML(V)$ of all partial isomorphisms of the vector space $V$ is
an inverse submonoid of $\II_V$. Moreover, $ ML(V)$ is factorizable if
and only if $V$ is finite dimensional. 
\end{lemma}

A system of subspaces $\BB$ for a subgroup $G$ of $ GL(V)$ is a special
case of a system of subsets for $G$ regarded as a subgroup of 
of $\goth{S}_V$, so as above we can construct
 a factorizable inverse submonoid $M(G,\BB)$ of $ ML(V)$ with group of
units $G$ and idempotents, the partial identities $\ve_X^{}$ for $X\in \BB$.

On the other hand, if $F$ is a factorizable inverse submonoid of $ML(V)$,
then we know that $F = M(G,\BB)$ where $\BB = \{\, \dom\, \sigma \mid \sigma \in
M\}$; now the domain of every element in $F$ is a subspace of $V$, so
$\BB$ is, in fact, a system of subspaces.   

A {\em partial
reflection\/} of a vector space $V$ is defined to be the restriction of
a reflection  $s\in GL(V)$ to a subspace $X$ of $V$. We denote this
partial reflection by $s_X^{}$. A \textit{reflection monoid} is defined to
be a
factorizable inverse submonoid of $ ML(V)$ generated by partial reflections.

It is easy to see that the non-units in a reflection monoid $M
\subset ML(V)$ form a subsemigroup, and hence every unit of $M$ must be a
product of (full) reflections, that is, the group of units of $M$ is a
reflection group $W$. Indeed, if $S$ is the set of generating partial
reflections for $M$, let $S'\subset S$ be the subset
of full reflections. Then $W=\langle S'\rangle$. Also, since $M$ is
factorizable, it follows that $M = M(W,\BB)$ for a system of subspaces
for $W$.

If we choose a (subspace) system $\BB$ for a reflection group $W\subset  GL(V)$,
the units of $M(W,\BB)$ are generated by reflections. Any other element
has the form $\ve_X^{}g$ for some $X\in \BB$ and $g\in W$. Now $g =
s_1\dots s_k$ for some reflections $s_1,\dots,s_k$ and $\ve_X^{}s_1$ is
a partial reflection, so $\ve_X^{}g = (\ve_X^{}s_1)s_2\dots s_k$ is
 a product of
partial reflections. Thus $M(W,\BB)$ is a reflection monoid. 

Most of the elementary properties of reflection monoids appear in the
above discussion. For emphasis, we list them in
the following result.

\begin{proposition}\label{reflection_monoids:result400}
Every reflection monoid $M \subset  ML(V)$ has the form $M(W,\BB)$
where $W$ is the reflection group of units  and 
$\BB = \{ \dom\,\sigma \mid \sigma \in M\}$. Conversely, if $W\subset  GL(V)$
is a non-trivial reflection group and $\BB$ is a system of subspaces for $W$, then
$M(W,\BB)$ is a reflection monoid with group of units $W$.

In  $M = M(W,\BB)$ we have:

\noindent$(1)$ $\BB = \{\,\dom\,\sigma \mid \sigma \in M\} = \{\,\im \sigma \mid \sigma \in
M\}$,

\noindent$(2)$ $E(M) = \{ \ve_X^{} \mid X\in \BB\}$, and 

\noindent$(3)$ the inverse of $g_X^{}$ is $(g^{-1})_{Xg}^{}.$

Finally, $M(W,\BB)$ is finite if and only if $W$ and $\BB$ are finite.
\end{proposition}

Recall that in any monoid $M$, Green's relation $\mathscr{R}$ is defined
by the rule that $a\mathscr{R} b$ if and only if $aM = bM$. The relation
$\mathscr{L}$ is the left-right dual of $\mathscr{R}$; we define 
$\mathscr{H} = \mathscr{R} \cap \mathscr{L}$ and $\mathscr{D} =
\mathscr{R}\ \vee \mathscr{L}$. In fact, by \cite[Proposition 2.1.3]{Howie95}, 
$\mathscr{D} = \mathscr{R} \circ \mathscr{L} = \mathscr{L} \circ \mathscr{R}$.
Finally, $a\JJ b$ if and only if $MaM =
MbM$. In an inverse monoid, $a\mathscr{R} b$ if and only if $aa^{-1} =
bb^{-1}$ and similarly, $a\mathscr{L}b$ if and only if $a^{-1}a =
b^{-1}b$. More information on Green's relations can be found in
\cite{Howie95,Lawson98}. 

\begin{proposition} \label{Green's}
Let $\rho,\sigma$ be elements of the reflection monoid $M = M(W,\BB)$ with
$\rho = g_X^{}$ and $\sigma = h_Y^{}$ where $g,h \in W$ and $X,Y\in \BB$. Then
\begin{description}
\item [$(1).$] $\rho \mathscr{R} \sigma$ if and only if $X = Y$;
\item [$(2).$] $\rho \mathscr{L} \sigma$ if and only if $Xg = Yh$;
\item [$(3).$] $\rho \mathscr{D} \sigma$ if and only if $Y\in XW$;
\item [$(4).$] if $\BB$ consists of finite dimensional spaces, then
$\mathscr{J} = \mathscr{D}$.
\end{description}
\end{proposition}
\begin{proof}
$(1)$ and $(2)$ follow from \cite[Proposition 2.4.2]{Howie95} and the
well known fact that in $ ML(V)$ we have $\rho \mathscr{R} \sigma$ if and only
if $\dom\, \rho = \dom\, \sigma$, and $\rho \mathscr{L} \sigma$ if and only if
$\im \rho = \im \sigma$.

If $\rho \mathscr{D} \sigma$, then $\rho \mathscr{R} \tau \mathscr{L} \sigma$
for some $\tau \in M$, and it follows from (1) and (2) that 
$Y\in XW$. On the other hand, if $Y\in XW$, say $Y = Xk$ where $k\in W$,
then $Yh = Xkh$ so that $\sigma \mathscr{L} (kh)_X^{}$ by (2), and 
$(kh)_X^{} \mathscr{R} \rho$ by (1), whence $\rho \mathscr{D} \sigma$.

Certainly, $\mathscr{D} \subset \mathscr{J}$. If $\rho \mathscr{J}
\sigma$, then $\rho = \alpha\sigma\beta$ and $\sigma = \gamma\rho\delta$ for some
$\alpha,\beta,\gamma,\delta \in M$. Comparing domains gives $X\subset
Ya$ and $Y\subset Xb$ for some $a,b\in W$. If the dimensions are finite,
we get $Y = Xb$ so that $\rho \mathscr{D} \sigma$ by (3).
\qed
\end{proof}

We remark that although we have stated this result for reflection
monoids, an entirely analogous result holds for factorizable monoids in
general. 

\parshape=20 0pt\hsize 0pt\hsize 0pt\hsize 0pt\hsize 0pt\hsize 
0pt\hsize 0pt\hsize 0pt\hsize 0pt\hsize 0pt\hsize 0pt\hsize 0pt\hsize 
0pt.75\hsize 0pt.75\hsize 0pt.75\hsize 0pt.75\hsize 0pt.75\hsize
0pt.75\hsize 0pt.75\hsize
0pt\hsize
The power of realising reflection monoids in the form
$M(W,\BB)$ will be seen in the next  section where we produce a
wealth of examples and calculate their orders.  For now, we use the 
idea to give an example of a non-fundamental
reflection monoid in which the restriction of the Munn
representation to the group of units is one-one. 
(Of course, we have seen that $\JJ_X$
is not fundamental, but in this case there are distinct units which 
are $\mu$-related.)
First,  note that if $M
= M(W,\BB)$ is any reflection monoid, and $\aa \in M$ has domain $X$, then
for any $Y \in \BB$ we have 
\begin{equation} \label{equation}
               \aa^{-1}\ve_Y^{}\aa = \ve_{(Y\,\cap\,X)\aa}^{}.
\end{equation} 
Now let $V=\mathbb{R}^2$ and 
$\Phi\subset V$ the root system shown (stably isomorphic to the crytallographic
$G_2\subset\R^3$ of Table \ref{table:roots2}) 
and for $W$ take the subgroup of $W(\Phi)$ generated by $\rho$ and
$\tau$ where $\rho$ is a rotation through $2\pi/3$ and $\tau$ is the
reflection in the $y$-axis. Thus $W \cong \goth{S}_3$. 
The $\R$-spans of these roots, together
with $V$ and $0$, form a system (of subspaces) for $W$. The $\mu$-class
of the identity $\ve_V^{}$ is a normal subgroup of $W$ and so to show that
$\mu$ is trivial on $W$, it is enough to show that $\rho$ and $\ve_V^{}$
are not $\mu$-related 
This is clear from \eqref{equation} using any of the six lines
for $Y$. 
\vadjust{\hfill\smash{\lower 10pt
\llap{
\begin{pspicture}(0,0)(3,3)
\rput(1.2,1.5){\BoxedEPSF{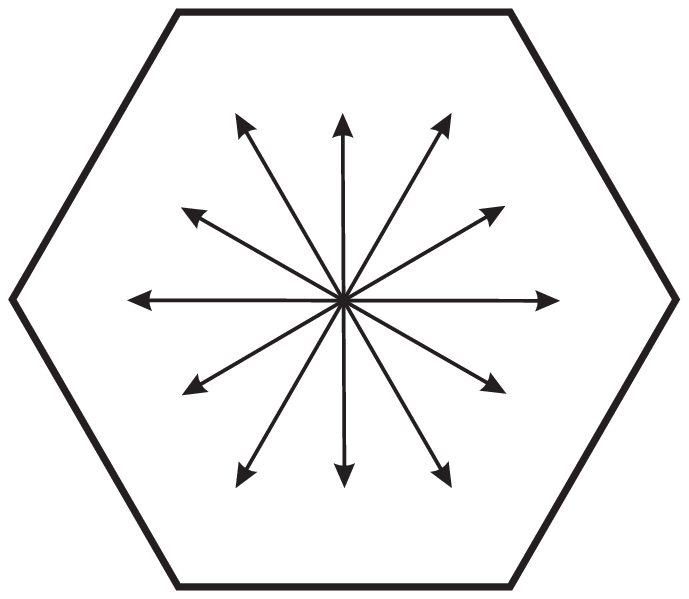 scaled 425}}
\end{pspicture}}}}\ignorespaces
On the other hand, letting 
$X$ be the $x$-axis, we see that $\tau_X$ and $\ve_X^{}$ are
distinct but $\mu$-related.

We now consider when two reflection monoids are
isomorphic. 
Let $W\subset GL(V),W'\subset GL(V')$ be reflection groups and $\BB,\BB'$
systems of subspaces for $W,W'$ respectively. We say that a vector space isomorphism 
$f:V\rightarrow V'$
induces an {\em isomorphism of reflection monoids\/} 
$f:M=M(W,\BB)\rightarrow M'=M(W',\BB')$ if 
$M'=f^{-1}Mf$.
It is easy to see that 
the map $\aa\mapsto f^{-1}\aa f$ is a monoid isomorphism $M\rightarrow M'$.

\begin{proposition}\label{reflection_monoids:result500}
$M$ and $M'$ are isomorphic reflection monoids if and only if
there is a vector space isomorphism  $f:V\rightarrow V'$ with $W'=f^{-1}W f$
and $\BB f=\BB'$, ie: $f:W\rightarrow W'$ is an isomorphism of reflection 
groups with $\BB f=\BB'$.
\end{proposition}

In particular, if the systems are the intersection lattices of hyperplane
arrangements, then
as $\rk X=\codim X$, an isomorphism of reflection monoids will induce
a bijection between the rank $k$ elements of the two systems.

\begin{proof}
If $f$ is an isomorphism of reflection monoids then the monoid isomorphism
$\aa\mapsto f^{-1}\aa f$ sends units to units, hence $W'\subset f^{-1}Wf$ with
$f^{-1}$ giving the reverse. If $X\in\BB$ then $\ve_X^{}\in M$, 
hence $Xf=f^{-1}\ve_X^{}f
\in M'$ giving $Xf\in\BB'$. Thus $\BB f\subset\BB'$ and $f^{-1}$
again gives the reverse. Conversely, if $f$ an isomorphism of the reflection
groups $W$ and $W'$ with $\BB f=\BB'$ and $\aa\in M$ then $\aa=g_X^{}$ for
$g\in W$ and $X\in\BB$ hence $f^{-1}\aa f=(f^{-1}gf)_{Xf}^{}$ with $Xf\in\BB'$
and $f^{-1}gf\in W'$ giving $f^{-1}\aa f\in M'$.
\qed
\end{proof}

We now proceed to find the orders of our reflection monoids, for which
the following result is straightforward but crucial.

\begin{theorem}\label{orders:result100}
Let $W\subset GL(V)$ be a reflection group and $\BB$ a system
for $W$. Then 
$$
|M(W,\BB)|=\sum_{X\in\BB}[W:W_X],
$$
where $W_X\subset W$ is the isotropy group of $X\in\BB$.
\end{theorem}

\begin{proof}
For $X\in\BB$ let $M(X)$ be the set of $\aa\in M(W,\BB)$ with
$\dom(\aa)=X$. Then $M(W,\BB)$ is the disjoint union of the 
$M(X)$ and so 
$|M(W,\BB)|=\sum_{X\in\BB}|M(X)|$. 
The 
elements of 
$M(X)$ are the partial isomorphisms obtained by restricting the elements
of $W$ to $X$, and $w_1,w_2\in W$ yield the same partial isomorphism if and only if they
lie in the same coset of the isotropy subgroup $W_X$. Thus, 
$|M(X)|=[W:W_X]$ and the result follows.
\qed
\end{proof}

Observe that the $M(X)$ of the proof is the $\mathscr{R}$-class containing the 
partial identity $\ve_X^{}$, so that the sum of Theorem \ref{orders:result100}
can be interpreted as a sum over $\mathscr{R}$-classes.
The proof also shows that the result is true for an arbitrary $G\subset GL(V)$, however,
 when $G$ is not a reflection group, it may not be so easy to calculate
the number of orbits and their sizes, and dealing with the isotropy
group $G_X$ may be difficult.
%
%
%

If $X,Y\in\BB$ lie in the same orbit of the $W$-action on $\BB$, then their
isotropy groups $W_X,W_Y$ are conjugate, and 
the sum in Theorem \ref{orders:result100} becomes 
\begin{equation}\label{orders:result300}
|M(W,\BB)|=|W|\sum_{X\in\,\Om} \frac{n_X^{} }{|W_X|},
\end{equation}
where $\Om$ is a set of orbit representatives for the $W$-action on $\BB$,
and $n_X^{}$ is the size of the orbit containing $X$.
Most of our applications of Theorem \ref{orders:result100} will use the form
(\ref{orders:result300}).


\section{Examples}\label{section:examples}

In this section we identify some important monoids
pre-existing in the literature as reflection monoids, and introduce some new
examples. In some cases the choices
are motivated by reflection groups that can be identified with other
common or garden variety groups.

\subsection{Boolean monoids}
\label{subsection:permutationmonoids}

We saw in \S\ref{section:reflectiongroups} that the classical Weyl 
groups have alternative
descriptions as groups of permutations, with $W(A_{n-1})\cong\goth{S}_n$,
$W(B_n)\cong\BBB_n=\goth{S}_n\ltimes\QQ(n)$ and 
$W(D_n)=\goth{S}_n\ltimes\QQ^+(n)$.

Much the same happens in the partial case. 
Let $W=W(\Phi)$ be a Weyl group as in Tables \ref{table:roots1}-\ref{table:roots2},
and $\BB=\langle\AA\rangle_W$ the Boolean system of 
\S\ref{subsection:permutationsystem}. Then the
the resulting reflection monoid $M(W,\BB)=M(\Phi,\BB)$ 
is called a {\em Boolean (reflection) monoid}.
Both $M(A_{n-1},\BB)$ and $M(B_n,\BB)$ can be identified with naturally occurring
permutation monoids.

Returning to the inverse monoids $\II_n$ and $\JJ_n$ of the previous
section, let $X = \{ 1,\dots,n\}$.
If $\{i,i+1\}\subset Y\subset X$,  let $\sigma_{i,Y}^{}$
be the partial permutation with domain and image $Y$, and whose effect on $Y$
is as the transposition $(i,i+1)$, ie: $\sigma_{i,Y}^{}$ interchanges $i$ and $i+1$,
fixes the remaining points of $Y$, and is undefined on $X\setminus Y$.
Similarly, let $\tau_{i,Y}^{}$ have domain and image $Y\cup-Y\subset X\cup -X$
with $i\in Y$ and effect $(i,-i)$ on $Y\cup -Y$; let $\mu_{i,Y}^{}$ 
have effect $(i,i+1)(-i,-(i+1))$ on
$Y\cup -Y$ for $\{ i,i+1\}\subset Y$. 

\begin{lemma}\label{examples:permutationmonoidsresult200}
Let $n\geq 3$.
$(1)$. The symmetric inverse monoid $\II_n$ is generated by the partial 
transpositions $\sigma_{i,Y}^{}$ for $1\leq i\leq n-1$ and $Y\subset X$.

\noindent$(2)$. The monoid of partial signed permutations $\JJ_n$ is generated by 
the $\tau_{i,Y}^{}$ and 
$\mu_{j,Y}^{}$ for $1\leq i\leq n$, $1\leq j\leq n-1$, $Y\subset X$.
\end{lemma}

\begin{proof}
For (1), we note that  $\goth{S}_n$ is generated by the full transpositions
$\sigma_{i,X}$, and that by (the proof of) \cite[Theorem 3.1]{Gomes87},
$\II_n$ is generated by any generating set for $\goth{S}_n$ together with any
partial permutation of rank $n-1$.

For (2), we recall that  $\JJ_n$ is factorizable so that every
element can be written as $\ve_{Y\cup -Y}\tau$ for some (full) signed
permutation $\tau$. Certainly $\BBB_n$ is generated by the $\tau_{i,X}$ 
  and $\mu_{i,X}$, so 
 it suffices to express
 $\ve_{Y\cup -Y}$ in terms of the proposed generating set. 
Writing $\ve_{i_1\ldots i_k}$ for $\ve_{Y\cup -Y}$ where 
$Y=X\setminus\{i_1,\ldots,i_k\}$,
we have $\ve_{i_1\ldots i_k}=\ve_{i_1}\ldots \ve_{i_k}$; hence it is enough to 
show that 
$\ve_{i}$ (for $1\leq i\leq n$) can be expressed in terms of the proposed generators.
  As $n\geq 3$, we have $\ve_n=\tau_{1,Y}^2$ where 
$Y=X\setminus\{n\}$,
and $\ve_j=\mu_{j,X}^{}\ve_{j+1}\mu_{j,X}^{}$ for $j<n$; hence
$\ve_1,\dots,\ve_n$ are generated by the $\tau_{i,X}$ 
  and $\mu_{i,X}$ as required. 
\qed
\end{proof}

Let $V$ be a Euclidean space with orthonormal basis $\{\xx_1,\ldots,\xx_n\}$,
and for $Y=\{i_1,\ldots,i_k\}\subset X=\{1,\ldots,n\}$, write 
$\langle Y\rangle$ for the span of $\{\xx_{i_1}\ldots,\xx_{i_k}\}$. If $\xx\in V$
let $s_{\xx}$ be the reflection in $\xx^\perp$ and $(s_\xx)_{\langle Y\rangle}$
the corresponding partial reflection.

\begin{proposition}\label{examples:permutationmonoids:result300}

$(1)$. The map 
$(s_{\xx_i-\xx_{i+1}})_{\langle Y\rangle}^{}\mapsto\sigma_{i,Y}^{}$
induces an isomorphism from the Boolean reflection monoid $M(A_{n-1},\BB)$ to the
symmetric inverse monoid $\II_n$.

\noindent$(2)$. The map 
$(s_{\xx_i})_{\langle Y\rangle}^{}\mapsto\tau_{i,Y}^{}$ and
$(s_{\xx_i-\xx_{i+1}})_{\langle Y\rangle}^{}\mapsto\mu_{i,Y}^{}$
induces an isomorphism from the Boolean reflection monoid $M(B_{n},\BB)$ to the
monoid of partial signed permutations $\JJ_n$.
\end{proposition}

\begin{proof}
As mentioned in \S \ref{section:reflectiongroups}, 
it is well known that when restricted to full
reflections, the map in (1) induces an isomorphism $\vphi:W(A_{n-1}) \to
\goth{S}_n$. For $g,h \in W(A_{n-1})$ and $Y \subset X$, it is clear that 
$g_{\langle Y \rangle}^{}=h_{\langle Y \rangle}^{}$ if and only if
$(g\vphi)_Y^{}= (h\vphi)_Y^{}$. Hence there is a bijection
$\ov{\vphi}:M(A_{n-1},\BB) \to \II_n$ extending $\vphi$ and given by 
$g_{\langle Y \rangle}^{}\ov{\vphi} = (g\vphi)_Y^{}$. It is easy to verify
that $\ov{\vphi}$ is an isomorphism which restricts to the map given in
(1). It follows that this map induces $\ov{\vphi}$ since the $\sigma_{i,Y}^{}$ 
generate $\II_n$. The proof of (2) is similar. 
\qed
\end{proof}

Unlike the Weyl group $W(D_n)$, there seems to be no nice interpretation of the 
reflection monoid $M(D_n,\BB)$ as a group of partial permutations. 
Now to the orders:

\begin{theorem}\label{orders:result400}
Let $\Phi_n$ be a root system of type $A_{n-1},B_n$ or $D_n$ as in Table
\ref{table:roots1} and $\BB$ the Boolean system for $W(\Phi_n)$.
Then the Boolean reflection monoids have orders,
$$
|M(\Phi_n,\BB)|=|W(\Phi_n)|\sum_{k=0}^n 
\left(\begin{array}{c}n\\k\end{array}\right) \frac{1}{|W(\Phi_k)|}.
$$
\end{theorem}

\begin{proof}
The $W$-action on $\BB$ is rank preserving and 
transitive on the rank $k$ elements,
with 
$\rk (X=\xx_{i_1}^\perp\cap\cdots\cap \xx_{i_k}^\perp)=\rk\{i_1,\ldots,i_k\}$ $=k$
(see \S\ref{subsection:permutationsystem}).
Thus the 
$X=\xx_{1}^\perp\cap\cdots\cap \xx_{k}^\perp$ for $0\leq k\leq n$
are orbit representatives, with 
$n_X$ 
the number of $k$ element subsets of $I$, and $W_X$ generated by the 
reflections $s_\vv$ for 
$\vv\in\Phi_n\cap X^\perp\cong\Phi_k$. The result now follows from (\ref{orders:result300}).
%
\qed
\end{proof}

By the conventions of \S\ref{section:reflectiongroups}
we have $|W(A_k)|=(k+1)!$, $|W(B_k)|=2^k k!$, $|W(D_0)|=1$,
and $|W(D_k)|=2^{k-1} k!$ for $k>1$, thus giving,
$$
\begin{tabular}{c|ccc}
\hline
$\Phi_n$&$A_{n-1}$&$B_n$&$D_n$\\\hline
&&&\\
$|M(\Phi_n,\BB)|$
&\vrule width 2 mm height 0 mm depth 0 mm
${\ds \sum_{k=0}^{n}\left(\begin{array}{c}n\\k\end{array}\right)^2
k!}$
&\vrule width 2 mm height 0 mm depth 0 mm
${\ds \sum_{k=0}^{n}2^k\left(\begin{array}{c}n\\k\end{array}\right)^2
k!}$
&\vrule width 3 mm height 0 mm depth 0 mm
${\ds 2^{n-1}n!+\sum_{k=1}^{n}2^k\left(\begin{array}{c}n\\k\end{array}\right)^2
k!}$\\
&&&\\\hline
\end{tabular}
$$
Notice that the given orders gel with the isomorphisms 
$M(A_{n-1},\BB)\cong\II_n$ and $M(B_n,\BB)\cong\JJ_n$ of Proposition
\ref{examples:permutationmonoids:result300} and the well known order 
of $\II_n$ (see eg: \cite[Chapter 5, Exercise 3]{Howie95}): 
one can independently choose a domain 
and image of size
$k$ for a partial permutation $\sigma\in\II_n$, with
there then being $k!$ partial permutations having the given domain and image; 
similarly for $\JJ_n$, there being $2^k k!$ partial signed permutations
with a given domain and image. One can also show, by thinking in terms
of partial signed
permutations, that the non-units of $M(B_n,\BB)$ and $M(D_n,\BB)$
coincide, which is why the orders of these reflection monoids are identical
except for the $k=0$ terms (recall that $\varnothing\subset I$ corresponds
to the ambient space $V\in\BB$). 


\subsection{The Renner monoids}\label{subsection:rennermonoids}

The theory of linear algebraic monoids was developed independently, and then
subsequently collaboratively, by Mohan Putcha and Lex Renner during the 1980's.
Among the chief achievements of the theory is the classification 
\cite{Renner85a,Renner85b} 
of the reductive monoids, and the formulation of a Bruhat decomposition
\cite{Renner86} for a reductive algebraic monoid, with the role of the 
Weyl group being played
by a certain finite factorizable inverse monoid, coined the \emph{Renner monoid\/}
by Solomon \cite{Solomon90}.

Thus the Renner monoids play the same role for algebraic monoids that the Weyl 
groups
play for algebraic groups, and in this section we investigate
to what extent the analogy 
continues further. 
Standard references on algebraic groups are 
\cite{Borel91,Humphreys75,Springer98}, and on algebraic monoids, the books of
Putcha and Renner \cite{Putcha88,Renner05}. We particularly recommend the
excellent survey of Solomon \cite{Solomon95}.

Throughout, $\F$ is an algebraically closed field. 
An \emph{affine\/} (or \emph{linear\/}) \emph{algebraic monoid\/} 
$\M$ over $\F$ is an affine algebraic variety 
together
with a morphism $\varphi:\M\times\M\rightarrow\M$ of varieties,
such that the product $xy=\varphi(x,y)$ gives $\M$ the structure of a monoid
(ie: $\varphi$ is an associative morphism of varieties and there is
a two-sided unit $1\in\M$ for $\varphi$).
We will assume that the monoid $\M$ is \emph{connected\/}, that is, the
underlying variety is irreducible, in which case the group $\G$ of units 
is a connected algebraic group with $\overline{\G}=\M$ (Zariski closure).
Adjectives normally applied to $\G$ are then transferred to $\M$; thus we have
semisimple monoids, reductive monoids, simply connected monoids, and so on

From now on, let $\M$ be reductive.
The key players, just as they are for algebraic groups, 
are the maximal tori $T\subset\G$ and their closures
$\overline{T}\subset\M$. Let $\XXX(T)$ be the character group of all
morphisms of algebraic groups $\chi:T\rightarrow\G_m$ (with 
$\G_m$ the multiplicative group of $\F$) and $\XXX(\overline{T})$
similarly the commutative monoid of morphisms of $\ov{T}$. Then $\XXX(T)$
is a free $\Z$-module, and restriction (together with the denseness
of $T$ in $\ov{T}$) embeds $\XXX(\ov{T})\hookrightarrow\XXX(T)$.

The Weyl group $W_\G=N_\G(T)/T$ of automorphisms of $T$ acts faithfully on 
$\XXX(T)$ via $\chi^g(t)=\chi(g^{-1}tg)$, thus realizing an injection
$W_\G\hookrightarrow GL(V)$ for $V=\XXX(T)\otimes\R$.
We will write $W$ for both the Weyl group and its image in $ GL(V)$. 
The non-zero weights $\Phi:=\Phi(\G,T)$ of the adjoint representation
$\G\rightarrow GL(\goth{g})$ form a root system with 
the Weyl group $W$ generated by reflections $s_\aa$ for
$\aa\in\Phi$ (with respect to a $W$-invariant bilinear form).

The Renner
monoid \cite{Renner86} $R_\M$ of $\M$ is defined to be
$R_\M=\overline{N_\G(T)}/T$, which turns out (although this is not obvious) to be
$N_\M(T)/T$, where $N_\M=\{x\in\M\,|\,xT=Tx\}$.
Just as $\II_n$ is the archetypal inverse monoid, and as 
$M(A_{n-1},\text{Boolean})$ it is the archetypal reflection monoid,
so in its incarnation as the rook monoid it is the standard example of a Renner
monoid, namely for $\M=\m_n(\F)$, the algebraic monoid of $n\times n$
matrices over $\F$. These monoids have been explicity described
in some other cases, for example, when $\M$ 
is the ``symplectic monoid'' $\text{MSp}_n(\F)=\ov{\F^*\sp_n(\F)}\subset\m_n(\F)$
\cite{Renner03}.


Suppose now $\M$ has a zero, and let $E=E(\overline{T})$ be the lattice of 
idempotents of $\overline{T}$, for $T$ a maximal torus in $\G$. Then
by the results of \cite[Chapter 6]{Putcha88}, $E$ is a graded lattice
with $\hat{0}$ and $\hat{1}$. Moreover, by \cite[Theorem 10.7]{Putcha88}, 
the Weyl group $W$ is the automorphism group of $E$, 
via $e^g=g^{-1}e g$,
and
by \cite[Remark 11.3(i)]{Putcha88}
the Renner monoid
$R_\M$ is the submonoid $\langle E,W\rangle\subset\TT_E$, of the
Munn semigroup $\TT_E$ of $E$. 
Thus by Proposition \ref{reflection_monoids:result50}, the Renner monoid 
has the form $M(W,\CC)$ where 
$\CC=\{Ex\,|\,x\in E\}$ is a system of subsets in $E$.

Before proceeding we summarize some basic facts about cones
from \cite[\S 1.2]{Fulton93}. If $V$ is a real space and $v_1,\ldots,v_s$
a finite set of vectors, then the convex polyhedral cone with generators $\{v_i\}$
is the set $\sigma=\sum\lambda_iv_1$ where $\lambda_i\geq 0$. The dual cone 
$\sigma^\vee\subset V^*$ consists of those $u\in V^*$ taking non-negative values
on $\sigma$. A face $\tau\subset\sigma$ is the intersection with $\sigma$ of the kernel $u^\perp$
of a $u\in\sigma^\vee$, and the faces form a meet semilattice $\FF(\sigma)$ under inclusion.
If $\tau\in\FF(\sigma)$, let $\ov{\tau}$ be the $\R$-span in $V$ of $\tau$, so that
if $\tau=\sigma\cap u^\perp$ for $u\in\sigma^\vee$, then $\sigma\cap\ov{\tau}=\tau$.
In particular, if 
$\bigcap\ov{\tau}_j=\bigcap\ov{\mu}_j$ in $V$ then 
$\bigcap\tau_j=\bigcap\mu_j$ in $\FF(\sigma)$.
Note that if $\{\tau_j\}\in\FF(\sigma)$ are faces of $\sigma$ then we have
$\ov{\tau}\subset\bigcap\ov{\tau}_j$ for
$\tau=\bigcap\tau_j$.

A cone is simplicial if it has a set $A=\{v_i\}$ of linearly 
independent generators. 
If $\tau_i$ is the cone on $\{v_1,\ldots,\wh{v}_i,\ldots,v_s\}$,
then $\tau_i=\sigma\cap u_i^\perp$, where $u_i$ is the vector corresponding to $v_i$ in the
dual basis for $V^*$. Thus $\tau_i$ is a face of $\sigma$, and the face lattice
$\FF(\sigma)$ is isomorphic to the Boolean lattice on the $1$-dimensional
faces $\R^+\cdot v_i$ of $\sigma$.
If
$\tau\in\FF(\sigma)$ corresponds to $A_\tau\subset A$ then 
$\tau_1\cap\tau_2$ corresponds to $A_{\tau_1}\cap A_{\tau_2}$, and
$\ov{\tau}=\R$-span of $A_\tau$.
In particular,
$\R\text{-span}\{\bigcap A_{\tau_j}\}
=\bigcap\{\R\text{-span}\,A_{\tau_j}\}$,
and so we have $\ov{\tau}=\bigcap\ov{\tau}_j$ 
when $\tau=\bigcap\tau_j$ for $\sigma$ simplicial.
Finally, a cone is strongly convex if the dual $\sigma^\vee$ spans $V^*$. Simplicial
cones are strongly convex. On the other hand, if $\dim V=2$, then
any strongly convex cone is simplicial \cite[1.2.13]{Fulton93}.

\begin{figure}
\begin{pspicture}(0,0)(15,4)
\rput(3,1){
\rput(.2,1.8){$\FF(\sigma)$}
\rput(1.8,1.8){$\FF(\sigma)$}
\rput(.2,.2){$E$}
\rput(1.8,.2){$E$}
\psline[linewidth=.1mm]{->}(.65,1.8)(1.35,1.8)
\psline[linewidth=.1mm]{->}(.4,.2)(1.6,.2)
\psline[linewidth=.1mm]{->}(.2,1.5)(.2,.5)
\psline[linewidth=.1mm]{->}(1.8,1.5)(1.8,.5)
}
\rput(7,-.5){
\rput(3.2,2.5){\BoxedEPSF{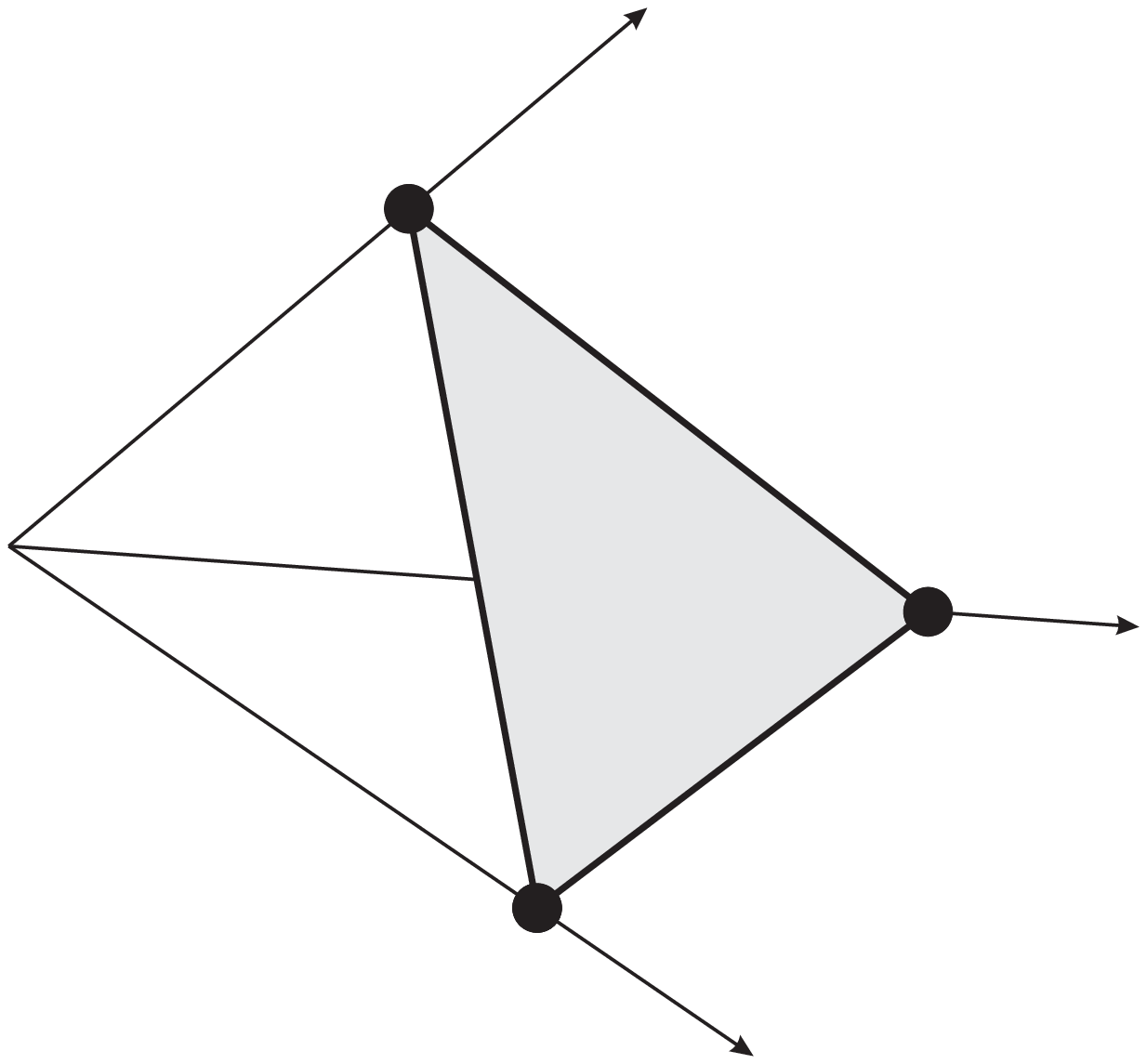 scaled 350}}
\rput(5,4){$\XXX(\ov{T})$}
\rput(3,0.8){${\scriptstyle 100}$}\rput(4.9,2.35){${\scriptstyle 010}$}
\rput(2.35,3.95){${\scriptstyle 001}$}\rput(4.15,1.55){${\scriptstyle 110}$}
\rput(3.95,3.1){${\scriptstyle 011}$}
\rput*(2.35,2.2){${\scriptstyle 101}$}\rput*(3.5,2.4){${\scriptstyle 111}$}
}
\end{pspicture}
\caption{(Left) The $W$-equivariant latice isomorphism $\FF(\sigma)\rightarrow E$:
the vertical arrows 
are the isomorphisms $\tau\mapsto e_\tau$, the top arrow is the map
$\tau\mapsto\tau g$ given by the $W$-action on the face lattice, and the
bottom arrow the map $e_\tau\mapsto e_\tau^g$ of the $W$-action on the idempotents.
(Right) the various ingredients for $\M=\m_3$, $\G= GL_3$,
$T=\d_3^*$, $\ov{T}=\d_3$, $\XXX(T)$ the free $\Z$-module on the 
characters 
$\chi_i(A)=A_{ii}$, 
$W_\G$ 
the permutation matrices,
and $E(\ov{T})$ the diagonal matrices with $0,1$-entries
(with the matrix $\text{diag}(a,b,c)$ represented by the string $abc$). 
The cone $\sigma$ is spanned by the shaded $2$-simplex (with vertices the basis
vectors $\chi_i$) with faces labelled by $E(\ov{T})$ to illustrate
the isomorphism $\FF(\sigma)\rightarrow E(\overline{T})$.
}
\label{section:examples:renner:figure50}
\end{figure}

Returning to algebraic monoids,
we may assume, by conjugating suitably, that the maximal
torus $T$ is a subgroup of the group
$\T_n$ of invertible diagonal matrices, where $n$ is the rank of $\G$. 
If $\chi_j$ is the restriction to $T$ of the
$j$-th coordinate function on $\T_n$, then the cone
$\sigma=\sum\R^+\chi_i\subset\XXX(T)\otimes\R$ is strongly convex. 
The dual cone $\sigma^\vee$ lives 
in the group of $1$-parameter subgroups of $T$.

This $\sigma$ has a number of nice properties.
Firstly, the character monoid $\XXX(\overline{T})=\sigma\cap\XXX(T)$.
Secondly, 
the Weyl group $W$, in its reflectional action on
$V$, acts on $\sigma$, and this induces an action 
$\tau\mapsto\tau g$
of $W$ on $\FF(\sigma)$.
Finally, the face lattice $\FF(\sigma)$ models the idempotents:
there is a lattice 
isomorphism
$\FF(\sigma)\rightarrow E(\ov{T})$, with $\tau\mapsto e_\tau$,
that is  $W$-equivariant 
with respect to the Weyl group actions, 
ie: for any $g\in W$, the diagram on the left 
of Figure \ref{section:examples:renner:figure50}
commutes. 
In short, $e_\tau^g=e_{\tau g}$ 
(Solomon \cite[Corollary 5.5]{Solomon95}, working with the dual
cone, has a lattice \emph{anti}-isomorphism $\FF(\sigma^\vee)\rightarrow E(\ov{T})$).

We can now define a reflection monoid using the Weyl group of
$\G$ and the convex polyhedral
cone $\sigma\subset\XXX(T)\otimes\Q$.
Let $\BB=\langle\ov{\tau}\,|\,\tau\in\FF(\sigma)\rangle_{W}$ be the system
for $W$ generated by the subspaces $\ov{\tau}$. As $W$ acts on the face 
lattice $\FF(\sigma)$, 
each $X\in\BB$ has the form $X=\bigcap\ov{\tau}_j$ for $\tau_j\in\FF(\sigma)$.
Call $M(W,\BB)$ the \emph{reflection monoid associated
to $\M$\/}.

Figure \ref{section:examples:renner:figure50}
depicts the situation for $\M=\m_3$.
The system $\BB$ is just the Boolean one generated by the coordinate hyperplanes
$\chi_i^\perp$, and the reflection monoid $M(W,\BB)$ is the symmetric inverse monoid
on the vertices of the $2$-simplex (hence, in this case, isomorphic to the Renner
monoid $R_\M$).

If $X=\bigcap\ov{\tau}_j\in\BB$, 
then the idempotents of $M(W,\BB)$ are products $\ve_X^{}=\prod\ve_j$ where
$\ve_j$ is the partial identity on $\ov{\tau}_j$, hence
any element of the reflection monoid has the form $\ve g=\prod\ve_j\cdot g$ for 
$g\in W$. Define a mapping $f:M(W,\BB)\rightarrow M(W,\CC)=R_\M$ by
$f(\ve g)=\prod e_j\cdot g$, where $e_j:=e_{\tau_j}\in E$. 

\begin{theorem}\label{section:examples:renner:result100}
Let $\M$ be connected reductive with $0$, 
$R_\M$ its Renner monoid, and $M(W,\BB)$ the associated reflection
monoid.
Then $f:M(W,\BB)\rightarrow R_\M$ is a surjective homomorphism, 
which is injective if and only if $\sigma\subset\XXX(T)\otimes\Q$ is a simplicial cone.
\end{theorem}

\begin{proof}
Let $X=\bigcap\ov{\tau}_j$, $Y=\bigcap\ov{\mu}_j$ and
$\ve_X^{}g_1=\ve_Y^{}g_2$ in the reflection monoid. Then 
$X=Y$ and $g_2^{}g_1^{-1}$ is in the isotropy group $W_X$ of $X$.
By intersecting the expressions for $X$ and $Y$ with $\sigma$ we
get $\bigcap\tau_j=\bigcap\mu_j$ in $\FF(\sigma)$, and so
$\prod e_{\tau_j}=\prod e_{\mu_j}$ in $E$, as these are the images
under the lattice isomorphism $\FF(\sigma)\cong E$.

Writing $e_j:=e_{\tau_j}$ and $\tau=\bigcap\tau_j$ from now on, 
it suffices, for $f$ to be well defined, to show that 
the elements $\prod e_j\cdot g_i$, $(i=1,2)$, give the same partial permutations
in $\II_E$, and this follows if $g_2^{}g_1^{-1}$ fixes the ideal $E(\prod e_j)$ pointwise.
Let $e_\kappa$ be in this ideal for some $\kappa\in\FF(\sigma)$, so that 
$\kappa\subset\tau$ by the isomorphism $\FF(\sigma)\cong E$, hence 
$\kappa\subset\ov{\kappa}\subset\ov{\tau}\subset\bigcap\ov{\tau}_j=X$. 
Thus, as $g_2^{}g_1^{-1}$
fixes $X$ pointwise, it fixes $\tau$ pointwise, giving
$$
e_\tau^{g_2^{}g_1^{-1}}=e_{\tau g_2^{}g_1^{-1}}=e_\tau,
$$
as the isomorphism $\FF(\sigma)\cong E$ is $W$-equivariant. Thus $f$ is
well defined. To see that it is a homomorphism, observe that
$\ve_X^{g^{-1}}=\ve_{Xg^{-1}}^{}$, where
$Xg^{-1}=\bigcap(\ov{\tau}g^{-1})$ and $X$ is as above. 
If $\tau_j\mapsto e_j$ via 
$\FF(\sigma)\cong E$, then $\tau_j g\mapsto e_{\tau_j g}=e_{\tau_j}^g$,
and so 
$$
\ve_X^{g^{-1}}\mapsto\prod e_{\tau_j}^{g^{-1}}=
\biggl(\prod e_{\tau_j}\biggr)^{g^{-1}}
$$
under $f$.
We then have
$$
\ve_X^{}g_1\cdot \ve_Y^{}g_2=\ve_X^{}\ve_Y^{g_1^{-1}}\kern-1mm g_1g_2
\stackrel{f}{\mapsto}
\prod e_{\tau_j}\biggl(\prod e_{\mu_j}\biggr)^{g_1^{-1}}\kern-2mm g_1g_2=
\prod e_{\tau_j}\cdot g_1\cdot \prod e_{\mu_j}\cdot g_2.
$$
Surjectivity is clear.

For the second part of the Theorem, let $\ov{\sigma}$ be the $\R$-span in $V$
of $\sigma$, where we must have $\dim\ov{\sigma}>2$ if $\sigma$ is not
simplicial.
There are then maximal faces $\tau_1,\tau_2\in\FF(\sigma)$
with $\tau_1\cap\tau_2=\{0\}$. As the $\ov{\tau}_i$ are hyperplanes
in $\ov{\sigma}$, the intersection
$\ov{\tau}_1\cap\ov{\tau}_2$ has codimension $2$ in $\ov{\sigma}$, hence
is non-zero.
If $\ve_i$ is the partial identity on $\ov{\tau}_i$ and $e_i:=e_{\tau_i}$, then
this translates into $\ve_1\ve_2\not=0$ in $M(W,\BB)$, but
$e_1e_2=0$ in $E$. In particular, the injectivity of $f$ fails, even on the idempotents.

On the other hand, if $\sigma$ is simplicial, let 
$\ve_X^{}g_1\mapsto \prod e_{\tau_j}\cdot g_1$, 
$\ve_Y^{}g_2\mapsto \prod e_{\mu_j}\cdot g_2$ with
$\prod e_{\tau_j}\cdot g_1=\prod e_{\mu_j}\cdot g_2$. As elements of $\II_E$ we have
$\prod e_{\tau_j}=\prod e_{\mu_j}$ and $g_1^{-1}g_2$ fixing the ideal
$E(\prod e_{\tau_j})$ pointwise. The lattice isomorphism then gives
$\tau=\bigcap\tau_j=\bigcap\mu_j=\mu$ and thus
$X=\bigcap\ov{\tau}_j=\ov{\tau}=\ov{\mu}=
\bigcap\ov{\mu}_j=Y$. 
If $\bigcap\tau_j$ is generated by the (independent) vectors
$v_1,\ldots,v_t$ and $\nu_i=\R^+\cdot v_i$,
then $e_{\nu_i}\in E(\prod e_{\tau_j})$ and so fixed by $g_1^{-1}g_2$.
Thus $\nu_i$ is also fixed, hence $X$ too, as it is spanned by such $\nu_i$.
Thus 
$g_1\ve_X^{}=g_2\ve_Y^{}$, and $f$ is injective.
\qed
\end{proof}

\begin{figure} 
\begin{pspicture}(0,0)(15,4)
\rput(-.5,0){
\rput(5.25,2){\BoxedEPSF{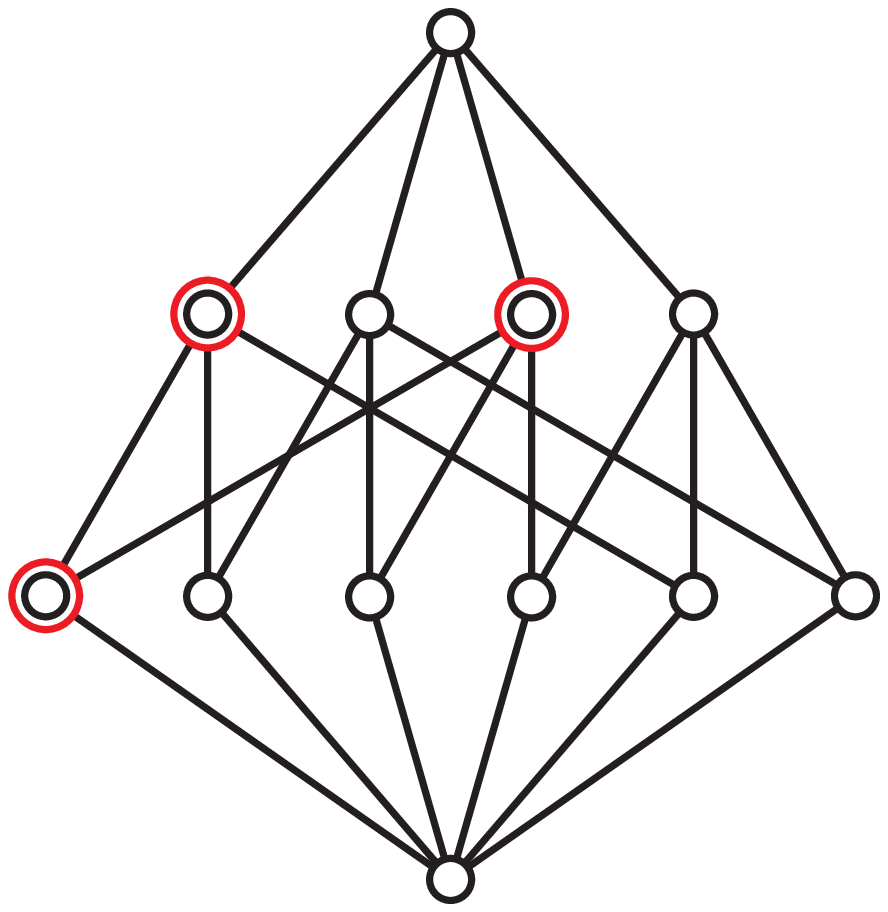 scaled 400}}
\rput(4,3.4){$\BB$}
\rput(3.9,2.6){{\red $\ov{\tau}_1$}}\rput(5.95,2.6){{\red $\ov{\tau}_2$}}
\rput(2.9,1.4){{\red $\ov{\tau}_1\cap\ov{\tau}_2$}}
}
\rput(7.35,2.2){$f$}
\psline[linewidth=.3mm]{->}(6.5,2)(8.2,2)
\rput(0.5,0){
\rput(10.6,3.4){$E(R_\M)$}
\rput(9.25,2){\BoxedEPSF{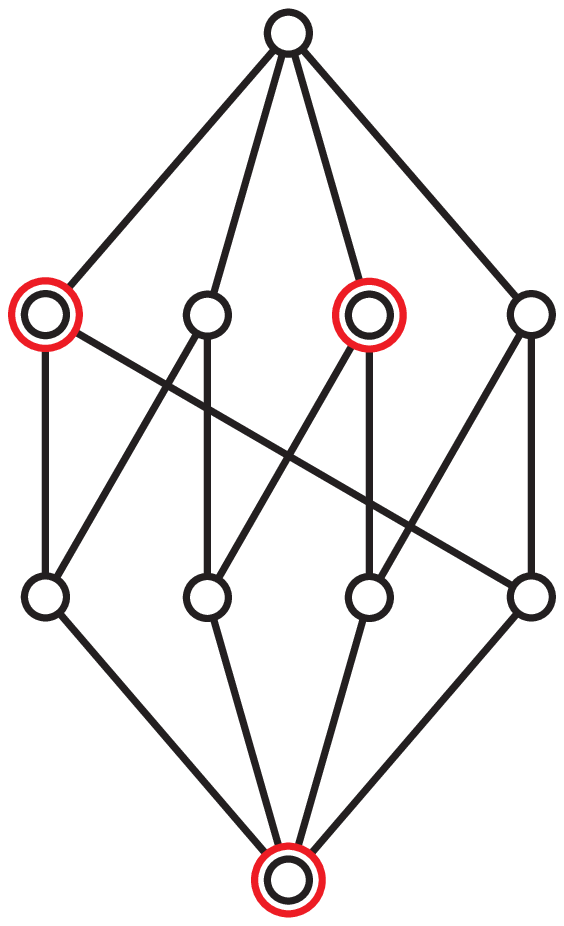 scaled 400}}
\rput(7.95,2.6){{\red $e_1$}}\rput(9.95,2.6){{\red $e_2$}}
\rput(9.95,.3){{\red $e_1\wedge e_2$}}
}
\end{pspicture}
\caption{The homomorphism $f$ of Theorem \ref{section:examples:renner:result100} 
need not be injective: if $\M=\ov{\text{Ad}(\G)\F^*}$ with $\G$ the adjoint simple group
of type $B_2$, then the lattice of idempotents of the associated reflection
monoid (left) contains non-zero elements mapping via $f$ to zero in the
lattice of idempotents of the Renner monoid
(right).}\label{section:examples:renner:figure100}
\end{figure}

As an illustration of the phenomenon in the last part of the proof, 
let $\M$ be the (normalization of) $\ov{\text{Ad}(\G)\F^*}$
for $\G$ the adjoint simple group of type $B_2$. Then
\cite[Example 3.8.3]{Renner85a}, $\dim(\XXX(T)\otimes\R)=3$
with $\sigma$ a cone on a square 
(see \cite[Figure 6]{Renner85a}). If $\tau_i$, $(i=1,2)$ are 
the cones on opposite, non-intersecting faces of the square, then
$\tau_1\cap\tau_2=\{0\}$, whereas $\ov{\tau}_1\cap\ov{\tau}_2$ is 
a $1$-dimensional subspace. Figure \ref{section:examples:renner:figure100}
gives the lattice of idempotents of the reflection monoid associated to
$\M$ (left) with a pair a $\ve_1\ve_2\not=0$ marked, mapping via $f$
to $e_1\wedge e_2=0$ (right).

Not only does the above homomorphism fail to be injective in this case,
but we can also show quite easily that $R_\M$ cannot be isomorphic to a
reflection monoid. For, suppose that $R_\M \cong M(W,\BB)$ where $\BB$
is a system of subspaces of a Euclidean space $V$ on which $W$ acts as a
reflection group. Since $W$ must be isomorphic to the group of units of
$R_\M$, we have $W = W(B_2)$. Hence four of the elements of order 2 in $W$
must be reflections. Also, the lattice $\BB$ must be isomorphic to the
lattice shown on the right in
Figure~\ref{section:examples:renner:figure100}. Moreover, if the bottom
element of $\BB$ is a non-zero subspace, we can factor it out to obtain
a lattice of subspaces with bottom element $\{0\}$.

Reading from left to right, let the atoms and coatoms of $\BB$ be
$U_0,U_1,U_2,U_3$ and $X_0,X_1,X_2,$ $X_3$ respectively. The intersection
of any two $U_i$'s is zero, as is the intersection of $X_0$ and $X_2$.
Hence for any choice of non-zero vectors $\uu_i \in U_i$ ($i =
0,1,2,3$), the set $\{ \uu_0^{},\dots,\uu_3^{} \}$ is linearly
independent.

The group of units of $R_\M$ is the automorphism group of $E(R_\M)$
where the action is by conjugation. Hence $W$ acting by conjugation  on
$\{ \ve_Y^{} \mid Y \in \BB\}$ gives all automorphisms of $E(M(W,\BB))$
and since $\ve_{Yg}^{} = g^{-1}\ve_Y^{}g$ for all $Y\in \BB$ and $g\in
W$, the same is true of the induced action of $W$ on $\BB$.

Automorphisms of $\BB$ are determined by their effect on the atoms. Let
$g,g'\in W$ be such that their actions give rise to the automorphisms
determined by interchanging $U_0$ with
$U_3$ and $U_1$ with $U_2$, and interchanging $U_0$ with
$U_1$ and $U_2$ with $U_3$ respectively.  Choose
$\uu_i\in U_i$ for $i=0,1$; then $\uu_0^{}g\in U_3$ and $\uu_1^{}g \in U_2$, so that 
$\{\uu_0^{},\uu_1^{},\uu_0^{}g,\uu_1^{}g\}$ is a basis for the subspace
it spans, say $U$. It is readily verified that $-1$ is an eigenvalue of
$g|_U^{}$ of multiplicity 2, so that $-1$ cannot be a simple eigenvalue
of $g$ itself. Thus $g$ (which has order 2) is not a reflection. 
Similarly, $g'$ is not a reflection. This is a contradiction since
there is only one element of order 2 in $W$ which is not a reflection. 
\vspace{1ex}

Our last result in this subsection is a negative one of sorts: 
if an inverse monoid $M$ is to be a reflection
monoid then we must have an injective homomorphism $M\hookrightarrow ML(V)$
with the units of $M$ a reflection group in $V$.

\begin{proposition}
Let $\M$ be connected with $0$ and $R_\M$ its Renner monoid. 
If $\rho:R_\M\rightarrow ML(V)$ is
faithful with $\rho(W_\G)$ a reflection group 
acting essentially on $V$, then $W_\G$ is not
of $(-1)$-type.
\end{proposition}

Thus at least one of the irreducible components of $W_\G$ must be
$A_n (n>1)$, $D_n$ ($n$ odd) or $E_6$.

\begin{proof}
It follows immediately that $W=\rho(W_\G)$ is a finite reflection group
acting essentially on $V$. In particular, $\rho$ is equivalent to the
reflectional representation of a Coxeter system $(W_\G,S)$, and 
if $W_\G$ is of $(-1)$-type, there is a $g\not= 1\in W_\G$ with $\rho(g)=-1$ on $V$. By
(\ref{equation}), $\rho(g)$ is $\mu$-related to $1\in\rho(R_\M)$, with the
resulting reflection monoid not fundamental.
\qed
\end{proof}

We conclude the subsection by mentioning that several authors have
calculated the orders of certain Renner monoids. The most general
results (which include all earlier ones) are in \cite{ZhenhengLi06}.

\subsection{Reflection arrangement monoids}
\label{subsection:reflectionarrangementreals}

Let $W\subset GL(V)$ be a reflection group and $\HH=L(\AA)$ the intersection lattice
of the arrangement $\AA$ of the reflecting hyperplanes of $W$.
The resulting $M(W,\HH)$ is called the {\em (reflection) arrangement
monoid\/} of $\AA$.

If $W=W(\Phi)$ we write
write $M(\Phi,\HH)$ for the arrangement monoid. If $\Phi\subset V$
and $\Phi'\subset V'$ are essential, then a root system isomorphism 
$f:\Phi\rightarrow\Phi'$ induces an isomorphism of reflection monoids
$M(\Phi,\HH)\rightarrow M(\Phi',\HH')$ where $\HH,\HH'$ are the lattices
of the reflection arrangements arising from $\Phi$ and $\Phi'$.
Thus we may talk of the arrangement monoids of types $A,B,\ldots$ etc, without
reference to the particular choice of root system, although we will usually have
in mind the $\Phi$ of \S\ref{section:reflectiongroups}.

Scrutinising these $\Phi$, we see that in
types $B$ and $F$, the Boolean system $\BB$ is properly contained in the arrangement
system $\HH$, thus the Boolean monoid $M(\Phi,\BB)$ is a proper submonoid of the
arrangement monoid $M(\Phi,\HH)$ in these cases. 
On the other hand,
an isomorphism of reflection monoids $M(\Phi,\BB)\rightarrow M(\Phi,\HH)$
would induce, by Proposition \ref{reflection_monoids:result500},
a bijection between the rank $k$ subspaces of the Boolean and 
arrangement systems.
For classical $\Phi$, the number of such subspaces in the arrangement systems
are
$$
\begin{tabular}{ccc}
\hline
$A$&$B$&$D$\\\hline
&&\\
$S(n,k)$&
$\sum_{i=0}^n 2^{i-k} \left(\begin{array}{c}n\\i\end{array}\right)S(i,k)$&
$\sum_{i\not=n-1} 2^{i-k} \left(\begin{array}{c}n\\i\end{array}\right)S(i,k)$\\
&&\\\hline
\end{tabular}
$$
where $S(n,k)$ is a Stirling number of the
second kind. As these numbers in the Boolean case are the number of ways of choosing 
$k$ objects from $n$, 
there is no isomorphism of reflection monoids
between $M(\Phi,\BB)$ and $M(\Phi,\HH)$ for these $\Phi$.


We now proceed to compute their orders, which in contrast to the Boolean
guys, we can do in both the classical and exceptional cases.
Recall from 
\S\ref{subsection:systemintersectionposet} that a partition of $n$
is a sequence of non-negative integers $\lambda=(\lambda_1,\ldots,\lambda_p)$
with $\sum\lambda_i=n$ and $\lambda_i\geq\lambda_{i+1}\geq 1$, and
if $b_i>0$ is the number of $\lambda_i$ equal
to $i$, then $b_\lambda=b_1!b_2!\ldots (1!)^{b_1}(2!)^{b_2}\ldots$

\begin{theorem}\label{orders:result500}
The arrangement monoid $M(A_{n-1},\HH)$ has order,
$$
|M(A_{n-1},\HH)|=
(n!)^2\sum_{\lambda}\frac{1}{b_\lambda \lambda_1!\ldots\lambda_p!},
$$ 
the sum over all partitions $\lambda$ of $n$.
\end{theorem}

The denominator of the sum in Theorem \ref{orders:result500} is largest for
the partition $\lambda=(n)$, which contributes $1/(n!)^2$, hence the
not {\em a priori\/} obvious fact that the sum is an integer.

\begin{proof}
This is another application of (\ref{orders:result300}), with 
by Proposition \ref{systems:result400},
the partitions of $n$ the orbit representatives, $n_{X(\Lambda)}^{}=n!/b_\lambda$
for $\lambda=\|\Lambda\|$, and $|W|=n!$.
The $W(A_{n-1})\cong\goth{S}_n$ action on $\HH$ is given by
$X(\Lambda)g(\sigma)=X(\Lambda\sigma)$ for $\sigma\in\goth{S}_n$, hence
$W_{X(\Lambda)}\cong\goth{S}_{\lambda_1}\times\cdots\times\goth{S}_{\lambda_p}$.
\qed
\end{proof}

Proceeding now to the type $B$ case, 
let 
$0\leq m\leq n$ be integers,
$$
c_{mn}\stackrel{\text{def}}{=}\sum_{i=0}^{\min\{m,n-m\}}
\left(\begin{array}{c}m\\i\end{array}\right)
\left(\begin{array}{c}n-m\\i\end{array}\right),
$$
and $\delta_{mn}\stackrel{\text{def}}{=}m!(n-m)!c_{mn}$. 
The following is more general than we need, but 
may be of independent interest:

\begin{proposition}\label{orders:result550}
The isotropy group $W_X\subset W(B_n)$ of the subspace
$X=X(\Delta,\Gamma,\Lambda)\in\HH$
has order 
$$
2^{m+p}\,m!\prod_{i=1}^{p}\delta_{\mu_i\lambda_i},
$$
where $m=|\Delta|$,
$\|\Lambda\|=(\lambda_1,\ldots,\lambda_p)$ 
for $\Lambda=\{\Lambda_1,\ldots,\Lambda_p\}$,
and $\mu_i=|\Gamma\cap\Lambda_i|$.
\end{proposition}

\begin{proof}
An element $g(\sigma,T)\in W(B_n)$ stabilizes $X$ precisely when 
$\Delta\sigma=\Delta$, $\Lambda_i\sigma=\Lambda_i$ and if 
$\Gamma_i=\Gamma\cap\Lambda_i$ and $T_i=T\cap\Lambda_i$, then for
each $1\leq i\leq p$, we have 
$(T_i\vartriangle\Gamma_i)\sigma=\Gamma_i$ or $\Lambda_i\setminus\Gamma_i$
(see \S\ref{subsection:systemintersectionposet}).
We are thus free in the first instance 
to choose a pair of $T_\Delta=T\cap\Delta$ and $\sigma$ any bijection
$\Delta\rightarrow\Delta$ (of which there are $2^m\,m!$) and the proof is 
completed by showing that the number of pairs of a $T_i$ and $\sigma_i$
(which is $\sigma$ restricted to $\Lambda_i$) is $2\delta_{\mu_i\lambda_i}$.
To have $(T_i\vartriangle\Gamma_i)\sigma_i=\Gamma_i$, it is clearly necessary that
$T_i\vartriangle\Gamma_i$ and $\Gamma_i$ have the same cardinality and 
conversely, if this is so then $\sigma_i$ can be the extension of any bijection
$T_i\vartriangle\Gamma_i\rightarrow\Gamma_i$. The $T_i\subset\Lambda_i$ for 
which $|T_i\vartriangle\Gamma_i|=|\Gamma_i|$ are precisely those subsets that
can be partitioned into two equal sized pieces, one contained in $\Gamma_i$
and the other in $\Lambda_i\setminus\Gamma_i$. The number of such is 
$c_{\mu_i\lambda_i}$ and for each one there are $\mu_i!$ bijections 
$T_i\vartriangle\Gamma_i\rightarrow\Gamma_i$, each one in turn extendable to
$(\lambda_i-\mu_i)!$ bijections $\sigma_i:\Lambda_i\rightarrow\Lambda_i$. 

The other possibility is that 
$(T_i\vartriangle\Gamma_i)\sigma_i=\Lambda_i\setminus\Gamma_i$, and as
$(\Lambda_i\setminus T_i)\vartriangle\Gamma_i=
\Lambda_i\setminus(T_i\vartriangle\Gamma_i)$, the map 
$T_i\mapsto\Lambda_i\setminus T_i$ is a bijection from the set of $T_i$
with $|T_i\vartriangle\Gamma_i|=k$ to the set of $T_i$ with 
$|T_i\vartriangle\Gamma_i|=\lambda_i-k$. The result is that there are
$c_{\mu_i\lambda_i}$ subsets $T_i$ with 
$|T_i\vartriangle\Gamma_i|=|\Lambda_i\setminus\Gamma_i|$, and 
$(\lambda_i-\mu_i)!\mu_i!$ bijections $\sigma_i:\Lambda_i\rightarrow\Lambda_i$
extending bijections $T_i\vartriangle\Gamma_i\rightarrow\Lambda_i\setminus\Gamma_i$. 
\qed
\end{proof}

For a partition $\lambda=(\lambda_1,\ldots,\lambda_p)$, let
$d_\lambda=4^p\,b_\lambda\lambda_1!\ldots\lambda_p!$

\begin{theorem}\label{orders:result600}
The arrangement monoid $M(B_n,\HH)$ has order,
$$
|M(B_n,\HH)|=
2^{2n-1}(n!)^2\sum_{m,\lambda}
\frac{1}{4^{m}\,d_\lambda},
$$ 
the sum over all pairs $(m,\lambda)$ where $0\leq m\leq n$ is an integer
and $\lambda$ is a partition of $n-m$.
\end{theorem}

\begin{proof}
Observe by Proposition \ref{systems:result500} that the orbit of the subspace
$X(\Delta,\Gamma,\Lambda)$ is determined by $m=|\Delta|$ and 
the partition $\lambda=\|\Lambda\|$ of $n-m$, 
with $\Gamma$ playing no role. We thus choose
$\Gamma=\varnothing$ in each orbit, 
and apply (\ref{orders:result300}) to $X(\Delta,\varnothing,\Lambda)$,
with $|W|=2^n n!$,
$$
n_X^{}=2^{n-m-p}\left(\begin{array}{c}n\\n-m\end{array}\right)
\frac{(n-m)!}{b_\lambda}\text{ and }
|W_X|=2^{m+p}\,m!\prod_{i=1}^{p}\lambda_i!,
$$
the last by Proposition \ref{orders:result550}.
\qed
\end{proof}

For the arrangement monoid of type $D$, 
the intersection lattice $\HH$ of the arrangement of reflecting 
hyperplanes is a sublattice of the type $B$ one.
It then suffices to compare the isotropy
groups in $W(B_n)$ and $W(D_n)$ of an $X\in\HH$.

\begin{proposition}\label{orders:result700}
If $\HH$ is the intersection lattice of the reflection arrangement for 
$W(D_n)$ and $X=X(\Delta,\Gamma,\Lambda)\in\HH$, then
the isotropy groups $W_X\subset W(D_n),W'_X\subset W(B_n)$
coincide when $\Delta=\varnothing$
and each $\lambda_i$ is even, otherwise $W_X$ has index $2$ in $W'_X$.
\end{proposition}

\begin{proof}
The index of $W_X$ in $W'_X$ is at most $2$ as $W_X=W(D_n)\cap W'_X$ with 
$W(D_n)$ of index two in $W(B_n)$. Thus either $W_X$ has index $2$ in $W'_X$
or the isotropy groups coincide, with the latter happening precisely when 
$Xg(\sigma,T)=X$ for $g(\sigma,T)\in W(B_n)$ implies that 
$g(\sigma,T)\in W(D_n)$, ie: that $|T|$ is even.
It is easy to check that this happens if and only if $\Delta=\varnothing$
and each $\lambda_i$ is even.
\qed
\end{proof}

\begin{theorem}\label{orders:result800}
The arrangement monoid $M(D_n,\HH)$ has order,
$$
|M(D_n,\HH)|=4^{n-1}(n!)^2
\sum_{m,\lambda} \frac{\varepsilon_{m,\lambda}^{}}{4^{m}\,d_{\lambda}}
,
$$ 
the sum over all 
pairs $(m,\lambda)$ where $0\leq m\leq n$ is an integer $\not= 1$
and $\lambda=(\lambda_1,\ldots,\lambda_p)$ is a partition of $n-m$,
with $\varepsilon_{m,\lambda}^{}=1$ if $m=0$ and each $\lambda_i$ is even, and
$\varepsilon_{m,\lambda}^{}=2$ otherwise.
\end{theorem}

\begin{proof}
Apply Propositions \ref{systems:result600}, \ref{orders:result550} and 
\ref{orders:result700} to (\ref{orders:result300}).
\qed
\end{proof}

The orders of the arrangement monoids for the exceptional Weyl groups
are calculated directly from (\ref{orders:result300}) and the data in Table 
\ref{table:orbitdata}. 

\begin{proposition}\label{orders:result900}
The orders of the exceptional arrangement monoids are 
$$
\begin{tabular}{c|ccccc}
\hline
$\Phi$&$G_2$&$F_4$&$E_6$&$E_7$&$E_8$\\\hline
&&&&&\\
$|M(\Phi,\HH)|$
&$7^2$&$11\cdot 4931$&$2^4\cdot 5^2\cdot 40543$
&$3\cdot 113\cdot 24667553$
&$11\cdot 79\cdot 55099865069$\\
&&&&\\\hline
\end{tabular}
$$
\end{proposition}


{}

\end{document}